\documentclass[a4paper,11pt]{amsart}

\usepackage[margin=1in]{geometry}
\usepackage{amssymb,mathrsfs}
\usepackage{mathabx}

\makeatletter
\@namedef{subjclassname@2020}{%
  \textup{2020} Mathematics Subject Classification}

\theoremstyle{plain}
\newtheorem{Theorem}{Theorem}
\newtheorem{Proposition}{Proposition}[section]
\newtheorem{Lemma}[Proposition]{Lemma}
\newtheorem{Corollary}[Proposition]{Corollary}

\theoremstyle{remark}
\newtheorem{Remark}[Proposition]{Remark}
\newtheorem*{Notations}{Notations}

\makeatother\numberwithin{equation}{section}
\everymath{\displaystyle}
\allowdisplaybreaks

\newcommand{\C}{\mathbb{C}}
\newcommand{\N}{\mathbb{N}}
\newcommand{\R}{\mathbb{R}}
\newcommand{\defeq}{\overset{\mathrm{def}}{=}}
\newcommand{\e}{\varepsilon}
\newcommand{\wt}[1]{\widetilde{#1}}
\let\Re\undefined
\DeclareMathOperator{\Re}{Re}
\let\Im\undefined
\DeclareMathOperator{\Im}{Im}

\renewcommand{\SS}{\mathcal{S}}
\newcommand{\HH}{\mathcal{H}}
\newcommand{\KK}{\mathcal{K}}
\newcommand{\PP}{P}
\newcommand{\QQ}{Q}

\def\<#1\>{\left\langle #1\right\rangle}
\newcommand{\logestimate}{\lesssim \log\bigl(1+|\xi|\bigr)}

\begin{document}

\title{On some hyperbolic equations of~third order}

\author[F.~Colombini]{Ferruccio Colombini}
\address{F.~Colombini,
          Dipartimento di Matematica,
          Universit\`{a} di~Pisa,
          Largo Bruno Pontecorvo 5,
          56127 Pisa,
          Italy}
\email{ferruccio.colombini@unipi.it}

\author[T.~Gramchev]{Todor Gramchev$^\dagger$}

\author[N.~Orr\`{u}]{Nicola Orr\`{u}}
\address{N.~Orr\`{u},
          Liceo Scientifico A.~Pacinotti,
          Via Liguria 9,
          09121 Cagliari,
          Italy}
\email{orru.nicola@yahoo.it}

\author[G.~Taglialatela]{Giovanni Taglialatela}
\address{G.~Taglialatela,
          Dipartimento di Economia e Finanza,
          Universit\`{a} di~Bari ``Aldo Moro'',
          Largo Abbazia S. Scolastica,
          70124 Bari,
          Italy}
\email{giovanni.taglialatela@uniba.it}

\thanks{This paper was started before the demise of T.~Gramchev and carried out in his memory.}

\begin{abstract}
We give sufficient conditions for the well-posedness in $\mathcal{C}^\infty$
of the Cauchy problem for third order equations with time dependent coefficients.
\end{abstract}

\maketitle

\tableofcontents

\section{Introduction}

In this paper we study the Cauchy Problem in~$\mathcal{C}^\infty$
for some weakly hyperbolic equation of~third order.
We are interested in Levi conditions,
that is in conditions on lower order terms
which ensure the well-posedness of~the Cauchy Problem.

In the case of~strictly hyperbolic equations of~order~$m$
(that is when the characteristic roots are real and distinct,
$m$ is a natural number) Petrowski~\cite{Petrowski1937} proved well-posedness
of the Cauchy Problem in~$\mathcal{C}^\infty$ for any lower order term.
Then Oleinik~\cite{Oleinik1970} studied weakly hyperbolic equations of~second order
(that is the two characteristic roots are real but may coincide)
with $\mathcal{C}^\infty$ coefficients and lower order terms
and gave some sufficient conditions for well-posedness
of the Cauchy Problem.
Nishitani~\cite{Nishitani1980} found necessary and sufficient conditions
for second order equations when there is only one space variable
and the coefficients are analytic.
In the papers~\cite{CDS} and~\cite{CJS}
some second order hyperbolic equations with coefficients depending only on~$t$
(in many space variables) were studied.
They studied the Cauchy problem both in~$\mathcal{C}^\infty$ and in Gevrey classes.
They studied the case of~$\mathcal{C}^\infty$ and
of analytic coefficients and gave some sufficient conditions
on the lower order terms for the well-posedness of~the Cauchy problem
(we~call them~\emph{logarithmic conditions}).

We~generalize these conditions to equations of~third order
with time dependent non smooth coefficients.
Recently Wakabayashi~\cite{Wakabayashi2015} has studied this problem,
obtaining results similar to~ours.
He~considers operators with double characteristics
and operators of~third order.
In~this case his sub-principal symbol is the same as ours,
while the sub-sub-principal symbol (of~order~1) is different.
His conditions on lower order terms
are similar to ours in one space variable.
In~general his conditions imply our logarithmic conditions.
There are two important differences between our works:
Wakabayashi supposes that the coefficients of~the equations are analytic,
while we~suppose that they are~$\mathcal{C}^2$
and his conditions are pointwise while our conditions are~integral.
We don't know if they are equivalent in many space variables,
when the coefficients are analytic.
Finally our conditions can be expressed simply
in terms of~the coefficients of~the operator
(at least if the coefficients are analytic).

\medskip

We have the following results:

\begin{Theorem}[\cite{CJS}, \cite{DT}, \cite{DKS}] \label{T-1}
Let's consider a second order equation
\[
\partial_t^2 u
   +  \sum_{j=1}^d a_j(t)\partial_t\partial_{x_j} u
   +  \sum_{j,h=1}^d b_{jh}(t) \partial_{x_j}\partial_{x_h} u
   + c_0(t) \partial_t u + \sum_{j=1}^d c_j(t)\partial_{x_j} u + d(t) u = f \, .
\]

Suppose that the coefficients are real and $\mathcal{C}^\infty$ in~$t$
and don't depend on the space variables.

Suppose that the symbol of~the principal part
\[
L(t,\tau,\xi) = \tau^2 + a(t,\xi)\tau + b(t, \xi) \, ,
\]
where $a(t,\xi) = \sum_{j=1}^d a_j(t)\xi_j$,
$b(t,\xi) =  \sum_{j,h=1}^d b_{jh}(t)\xi_j\xi_h$,
has real zeros in~$\tau$ for any $\xi\in\R^d$,
$t\in [0,T]$:
\[
\Delta(t,\xi) = a^2(t,\xi) - 4 b(t,\xi) \ge 0
\]
\textup(weak hyperbolicity\textup).

If
\[
\int_0^T
   \frac{\bigl|\partial_t\Delta (t,\xi)\bigr|}
        {\Delta (t,\xi)+1}\, dt
  \le  C\log |\xi|
\]
for any $\xi\in\R^n$ with $|\xi| \ge C_1>1$
\textup(this condition is automatically satisfied if the coefficients
are~analytic\textup),
and if
\begin{equation} \label{E-54}
\int_0^T
  \frac{\Bigl|-1/2 c_0(t) \sum_{j=1}^n a_j(t) \xi_j
               + \sum_{j=1}^n c_j(t) \xi_j
               - 1/2 \sum_{j=1}^n \partial_t a_j(t) \xi_j\Bigr|}
       {\sqrt{\Delta (t,\xi) +1\,}}\, dt
  \le  C \log|\xi|
\end{equation}
for any $\xi$ with $|\xi| \ge C_1>1$ \textup(logarythmic condition\textup),
then the Cauchy problem is well posed in~$\mathcal{C}^\infty$.
\end{Theorem}

Now we consider a third order equation
\begin{equation} \label{E-third}
\partial_t^3 u
  +  \sum_{j=0}^2\sum_{|\alpha|\le 3-j} a_{j,\alpha}(t)\partial_t^j\partial_x^\alpha u
     = f \, ,
\end{equation}
with initial conditions
\begin{equation} \label{E-IC}
u(0,x) = u_0 \, , \quad
\partial_t u(0,x) = u_1 \, , \quad
\partial_t^2 u(0,x) = u_2 \, .
\end{equation}

Let
\begin{align*}
L(t,\tau,\xi)
&  \defeq  \tau^3 + \sum_{j+|\alpha|=3} a_{j,\alpha}(t)\tau^j\xi^\alpha \, ,  \\
M(t,\tau,\xi)
&  \defeq  \sum_{j+|\alpha|=2} a_{j,\alpha}(t) \tau^j\xi^\alpha \, ,  \\
N(t,\tau,\xi)
&  \defeq  \sum_{j+|\alpha|=1} a_{j,\alpha}(t) \tau^j\xi^\alpha \, ,  \\
p(t)
&  \defeq  a_{0,0}(t) \, ,
\end{align*}
so that equation~\eqref{E-third} can be rewritten as
\[
L(t,\partial_t,\partial_x)u + M(t,\partial_t,\partial_x)u
  +  N(t,\partial_t,\partial_x)u + p(t)u = f \, .
\]

We assume that the coefficients of~$L$ belong to $\mathcal{C}^2\bigl([0,T]\bigr)$,
those of~$M$ and~$N$ belong to $\mathcal{C}^1\bigl([0,T]\bigr)$,
whereas $p(t)$ belongs to $L^\infty\bigl([0,T]\bigr)$.

The principal part $L(t,\tau,\xi)$,
as a polynomial in~$\tau$, has only real roots:
\[
\tau_1(t,\xi)\le \tau_2(t,\xi)\le \tau_3(t,\xi)
\]
for any $t,\xi$, (weak hyperbolicity).
This is equivalent to say that the~\emph{discriminant of~$L$}
is~nonnegative:
\begin{align*}
\Delta_L(t,\xi)
&  \defeq  \bigl(\tau_1(t,\xi)-\tau_2(t,\xi)\bigr)^2
           \bigl(\tau_2(t,\xi)-\tau_3(t,\xi)\bigr)^2
           \bigl(\tau_3(t,\xi)-\tau_1(t,\xi)\bigr)^2  \\
&  =  A_1^2(t,\xi) A_2^2(t,\xi) - 4 A_2^3(t,\xi) - 4 A_1^3(t,\xi) A_3(t,\xi)  \\
&  \qquad    +  18 A_1(t,\xi) A_2(t,\xi) A_3(t,\xi) - 27 A_3^2(t,\xi) \ge 0 \, ,
\end{align*}
where
\begin{equation} \label{E-Aj}
A_j(t,\xi)
  \defeq  \sum_{|\alpha|=3-j} a_{j,\alpha}(t) \, \xi^\alpha \, .
\end{equation}

We set also
\begin{align*}
\Delta_L^{(1)}(t,\xi)
&  \defeq  \bigl(\tau_1(t,\xi)-\tau_2(t,\xi)\bigr)^2
           +  \bigl(\tau_2(t,\xi)-\tau_3(t,\xi)\bigr)^2
           +  \bigl(\tau_3(t,\xi)-\tau_1(t,\xi)\bigr)^2  \\
&  =  2\,[A_1^2(t,\xi) - 3A_2(t,\xi)] \, .
\end{align*}

Note that if
$\Delta_L(\overline{t},\overline{\xi}) = 0$
and
$\Delta_L^{(1)}(\overline{t},\overline{\xi}) \ne 0$
then $L$ has a double root
(for example $\tau_1(\overline{t},\overline{\xi})=\tau_2(\overline{t},\overline{\xi})$
and $\tau_1(\overline{t},\overline{\xi})\ne\tau_3(\overline{t},\overline{\xi})$)
whereas if
$\Delta_L^{(1)}(\overline{t},\overline{\xi}) = 0$
then $L$ has a triple root:
$\tau_1(\overline{t},\overline{\xi})
  =  \tau_2(\overline{t},\overline{\xi})
  =  \tau_3(\overline{t},\overline{\xi})$.

Note that (cf.~Lemma~\ref{L-A-1})
\[
\Delta_L^{(1)}(t,\xi)
  =  \frac{9}{2} \, \Delta_{\partial L}(t,\xi) \, ,
\]
where
\[
\Delta_{\partial L}(t,\xi)
  =  \bigl(\sigma_1(t,\xi) - \sigma_2(t,\xi)\bigr)^2
\]
is the discriminant of~the polynomial
\[
\partial_\tau L(t,\tau,\xi)
  \defeq  3\,\tau^2 + 2\,A_1(t,\xi) \, \tau + A_2(t,\xi)
     =    3\,\bigl(\tau-\sigma_1(t,\xi)\bigr) \, \bigl(\tau-\sigma_2(t,\xi)\bigr) \, .
\]

\begin{Notations}
In the following we note
\begin{align*}
\SS_2
&  =  \Bigl\{ \ (1,2) \ , \ (2,3) \ , \ (3,1) \ \Bigr\}  \\
\SS_3
&  =  \Bigl\{ \ (1,2,3) \ , \ (2,3,1) \ , \ (3,1,2) \ \Bigr\} \, .
\end{align*}

Let $f(t,\xi)$ and $g(t,\xi)$ be positive functions,
we~will write $f\lesssim g$ (or, equivalently $g\gtrsim f$)
to~mean that there exists a positive constant $C$ such that
\[
f(t,\xi)
  \le  C \, g(t,\xi) \, ,
\quad
\text{for any $(t,\xi)\in[0,T]\in\R^n$} \, .
\]

Similarly,
we~will write $f\approx g$ to~mean that
$f \lesssim g$ and $g\lesssim f$.
\end{Notations}

These notations will make the~formulas more readable
and will allow us to~focus only on the~important terms
of~the~estimates.

\smallskip

Consider the auxiliary polynomial
\begin{equation} \label{E-cL}
\mathcal{L}(t,\tau,\xi)
   \defeq  L(t,\tau,\xi) - \partial_t^2 L(t,\tau,\xi) \, ,
\end{equation}
we can prove (see Lemma~\ref{L-Jannelli} below)
that its roots $\lambda_j(t,\xi) $ are real and distinct for $\xi\ne0$,
and, there exist positive constants $C_1$ and $C_2$ such that
\begin{align*}
\bigl|\lambda_j(t,\xi) - \tau_j (t,\xi)\bigr|
&  \le  C_1 \, ,  \\
\bigl|\lambda_j(t,\xi) - \lambda_k(t,\xi)|
&  \ge  C_2 \, ,
\end{align*}
for all $(t,\xi) \in [0,T]\times\R^n\setminus\{0\}$
and~$(j,k)\in\SS_2$.

We denote by~$\mu_1$ and~$\mu_2$ the roots of~$\partial_\tau \mathcal{L}(t,\tau,\xi)$.

Define the symbols
\begin{align}
\widecheck{M}(t,\tau,\xi)
&  \defeq  M(t,\tau,\xi) - \frac{1}{2} \partial_t\partial_\tau L(t,\tau,\xi) \, , \label{E-cM} \\
\widecheck{N}(t,\tau,\xi)
&  \defeq  N(t,\tau,\xi) - \frac{1}{2} \partial_t\partial_\tau M(t,\tau,\xi)
             + \frac{1}{12} \partial_t^2\partial_\tau^2 L(t,\tau,\xi) \, . \label{E-cN}
\end{align}

We can now state our main result.

\begin{Theorem} \label{T-2}
Assume that
\begin{gather}
\int_0^T \sum_{(j,k)\in\SS_2}
  \frac{\bigl|\partial_t \lambda_j (t,\xi) - \partial_t \lambda_k (t,\xi)\bigr|}
        {\bigl|\lambda_j (t,\xi) - \lambda_k (t,\xi)\bigr|} \, dt
  \logestimate \, , \label{E-1}  \\
\int_0^T \sum_{(j,k)\in\SS_2}
  \frac{\bigl|\partial_t^2 \lambda_j (t,\xi) - \partial_t^2 \lambda_k (t,\xi)\bigr|}
       {\bigl|\partial_t\lambda_j(t,\xi) - \partial_t \lambda_k (t,\xi)\bigr|+1} \, dt
  \logestimate \, , \label{E-2}  \\
\int_0^T \sum_{j=1}^3
  \frac{\bigl|\partial_t \widecheck{M}\bigl(t,\lambda_j(t,\xi),\xi\bigr)\bigr|}
       {\bigl|\widecheck{M}\bigl(t,\lambda_j(t,\xi),\xi\bigr)\bigl|+1} \, dt
  \logestimate \, , \label{E-3}  \\
\int_0^T \sum_{j=1}^2
  \frac{\bigl|\partial_t \widecheck{N}\bigl(t,\mu_{j}(t,\xi),\xi\bigr)\bigr|}
        {\bigl|\widecheck{N}\bigl(t,\mu_{j}(t,\xi),\xi\bigr)\bigr|+1} \, dt
  \logestimate \, , \label{E-4}  \\
\int_0^T \sum_{(j,k,l)\in\SS_3}
  \frac{\bigl|\widecheck{M}\bigl(t,\lambda_j(t,\xi),\xi\bigr)\bigr|}
        {\bigl|\lambda_j(t,\xi) - \lambda_k(t,\xi)\bigr|
        \cdot \bigl|\lambda_j(t,\xi) - \lambda_l(t,\xi)\bigr|} \, dt
  \logestimate \, , \label{E-Mcl}  \\
\int_0^T
     \sum_{j=1}^2
     \sqrt{ \frac{\bigl|\widecheck{N}\bigl(t,\mu_{j}(t,\xi),\xi\bigr)\bigr|}
                 {\bigl|\mu_2(t,\xi) - \mu_1(t,\xi)\bigr|} \, } \, dt
  \logestimate \, , \label{E-Ncl}
\end{gather}

Then the Cauchy problem~\eqref{E-third}-\eqref{E-IC}
is well-posed in~$\mathcal{C}^\infty$.
\end{Theorem}

In the following, we will say that a function $f(t,\xi)$
verifies the~\emph{logarithmic condition} if
\[
\int_0^T \bigl| f(t,\xi) \bigr| \, dt
  \lesssim \log\bigl(1+|\xi|\bigr) \, .
\]

\begin{Remark}
Conditions~\eqref{E-1}, \eqref{E-2}, \eqref{E-3} and~\eqref{E-4}
are hypothesis on the regularity of~the coefficients.
Indeed if the coefficients are analytic
then they are satisfied, see~\textsection\ref{S-an}.

Conditions~\eqref{E-Mcl} and~\eqref{E-Ncl} are Levi conditions
on the lower order terms.
They are necessary if the coefficients of~the principal symbols are constant,
see~\textsection\ref{S-cc}.
\end{Remark}

\begin{Remark}
The hypothesis in Theorem~\ref{T-2}
can be expressed in terms of~the coefficients of~the operator.
This is possible either by expliciting the roots of~$\mathcal{L}$
and~$\partial_\tau\mathcal{L}$,
or by transforming Hypothesis \eqref{E-1}--\eqref{E-Ncl}
into symmetric rational functions of~the roots of~$\mathcal{L}$.
This will be developped in~\textsection\ref{S-EF}.
\end{Remark}

In the proofs, for the sake of~simplicity,
we~will omit the dependence on $t$ and~$\xi$
in the notations.

\medskip

The plan of~the paper is the following.
In~\textsection\ref{S-proof} we~will prove Theorem~\ref{T-2}.
In~\textsection\ref{S-EF}
we give some different forms of~the Levi conditions~\eqref{E-Mcl} and~\eqref{E-Ncl}. 
In~\textsection\ref{S-an}
we will show that if the coefficients
are analytic, then~\eqref{E-1}, \eqref{E-2}, \eqref{E-3} and~\eqref{E-4}
are satisfied.
In~\textsection\ref{S-anLevi}
we give some sufficient pointwise conditions
that are equivalent to ours in space dimension $n=1$.
Finally, in~\textsection\ref{S-cc}
we~show that the Levi conditions~\eqref{E-Mcl} and~\eqref{E-Ncl}
are necessary for the $\mathcal{C}^\infty$ well-posedness
if the coefficients of~the principal symbols are constant.

\section{Proof of~Theorem~\ref{T-2}} \label{S-proof}

\begin{Lemma}[\cite{Jannelli1989}] \label{L-Jannelli}
Consider the polynomial
\[
L_\e (t,\tau,\xi) = L(t,\tau,\xi) - \e^2 |\xi|^2\partial_\tau^2 L(t,\tau,\xi) \, .
\]

Its roots $\tau_{j,\e} (t,\xi)$ are real and distinct,
moreover
\begin{align*}
\bigl|\tau_{j,\e} -\tau_j\bigr|
&  \lesssim \e |\xi| \, ,  &
   &  j = 1,2,3 \, ,  \\
\bigl|\tau_{j,\e} - \tau_{k,\e}\bigr|
&  \gtrsim \e |\xi| \, ,  &
   &  (j,k)\in\SS_2 \, .
\end{align*}
\end{Lemma}

\begin{Remark}
By direct calculation (cf.~\cite[par.~1020]{JannelliT2014}),
the discriminant of~$L_\e$ is given by
\[
\Delta_{L_\e}
  =  \Delta_L + \frac{1}{2} \e^2 |\xi|^2 \, \Delta_{\partial_\tau L}^2
     + 36\,\e^4 |\xi|^4\,\Delta_{\partial_\tau L} + 864\,\e^6 |\xi|^6 \, ,
\]
where $\Delta_{\partial_\tau L}$ is the discriminant of~the polynomial $\partial_\tau L$.

Similarly
\begin{equation} \label{E-DDD}
\Delta_{\partial_\tau L_\e}
  =  \Delta_{\partial_\tau L} + 72\,\e^2 |\xi|^2 \, .
\end{equation}
\end{Remark}

We will take~$\e = 1/ |\xi|$.

We consider
\begin{align*}
L_{j,\e}(t,\tau,\xi)
&  =  \tau - i\tau_{j,\e}(t,\xi) \, ,  &
   &  j = 1,2,3 \, ,  \\
L_{jk,\e}(t,\tau,\xi)
&  =  \bigl(\tau - i\tau_{j,\e}(t,\xi)\bigr)
      \bigl(\tau - i\tau_{k,\e}(t,\xi)\bigr) \, ,  &
   &  (j,k)\in\SS_2 \, ,  \\
L_{123,\e}(t,\tau,\xi)
&  =  \bigl(\tau - i\tau_{1,\e}(t,\xi)\bigr)
      \bigl(\tau - i\tau_{2,\e}(t,\xi)\bigr)
      \bigl(\tau - i\tau_{3,\e}(t,\xi)\bigr)
   =  L_\e \, .
\end{align*}

We define also the operators
\begin{align}
\wt{L}_{jh,\e}(t,\partial_t,\xi)
&  =  \frac{1}{2}
      \bigl[ L_{j,\e}(t,\partial_t,\xi) \circ L_{h,\e}(t,\partial_t,\xi)
             +  L_{h,\e}(t,\partial_t,\xi) \circ L_{j,\e}(t,\partial_t,\xi) \bigr] \, ,
             \label{E-wtLL}
\intertext{%
for any $(j,h)\in\SS_2$,
and}
\wt{L}_{123,\e}(t,\partial_t,\xi)
&  =  \frac{1}{6}
      \sum_{\substack{j,h,l=1,2,3 \\ j\ne h\,,\,j\ne l\,,\,h\ne l}}
      L_{j,\e}(t,\partial_t,\xi) \circ
      L_{h,\e}(t,\partial_t,\xi) \circ
      L_{l,\e}(t,\partial_t,\xi) \, . \label{E-wtLLL}
\end{align}

If $\e=0$, we will write $L_{j}$, $L_{jh}$, $\wt{L}_{jh}$, \dots,
instead of~$L_{j,0}$, $L_{jh,0}$, $\wt{L}_{jh,0}$, \dots.

\begin{Lemma}
For any $(j,h)\in\SS_2$,
we have
\begin{equation} \label{E-tLjh}
\wt{L}_{jh,\e}
  =  L_{jh,\e}
      + \frac{1}{2} \, \partial_t \partial_\tau L_{jh,\e} \, ,
\end{equation}
and
\begin{equation} \label{E-23}
L_{j,\e} \circ L_{h,\e} - \wt{L}_{jh,\e}
  =  \frac{i}{2} \, (\tau_{j,\e}' - \tau_{h,\e}') \, .
\end{equation}
\end{Lemma}

\begin{proof}
As
\begin{align}
L_{j,\e} \circ L_{h,\e}
&  =  (\partial_t - i\tau_{j,\e})
      \circ  (\partial_t - i\tau_{h,\e})
   =  L_{jh,\e}
      - i \tau_{h,\e}' \, , \label{E-21}
\intertext{%
we have}
(\partial_t - i\tau_{j,\e})
&  \circ  (\partial_t - i\tau_{h,\e})
   +  (\partial_t - i\tau_{h,\e})
      \circ  (\partial_t - i\tau_{j,\e}) \notag  \\
&  =  2 L_{jh,\e}
      - i (\tau_{j,\e}' + \tau_{h,\e}') \label{E-22}  \\
&  =  2 L_{jh,\e}
      + (\partial_t\partial_\tau L_{jh,\e}) \, , \notag
\end{align}
from which~\eqref{E-tLjh} follows.

Identity~\eqref{E-23} follows from~\eqref{E-21} and~\eqref{E-22}.
\end{proof}

Note that
\begin{equation} \label{E-32}
L_{j,\e} v - L_{h,\e} v
  =  -i(\tau_{j,\e} - \tau_{h,\e})v \, ,
\end{equation}
whereas,
{}from~\eqref{E-22},
\begin{equation} \label{E-31}
\begin{split}
\wt{L}_{jh,\e}v-\wt{L}_{jl,\e}v
&  =  L_{jh,\e}v-L_{jl,\e}v
      - \frac{i}{2} (\tau_{j,\e}' + \tau_{h,\e}')v
      + \frac{i}{2} (\tau_{j,\e}' + \tau_{l,\e}')v  \\
&  =  -i(\tau_{h,\e}-\tau_{l,\e}) \, L_{j,\e}v
      - \frac{i}{2} (\tau_{h,\e}' - \tau_{l,\e}')v \, .
\end{split}
\end{equation}

\begin{Lemma}
We have
\begin{align}
\wt{L}_{123,\e}
&  =  L_{123,\e}
      +  \frac{1}{2} \, \partial_t \partial_\tau L_{123,\e}
      +  \frac{1}{6} \, \partial_t^2 \partial_\tau^2 L_{123,\e}
   \label{E-tL123}  \\
L_{1,\e} \circ L_{2,\e} \circ L_{3,\e} - \wt{L}_{123,\e}
&  =  \frac{i}{2} (\tau_1'-\tau_2') \, L_{3,\e}
      +  \frac{i}{2} (\tau_2'-\tau_3') \, L_{1,\e}
      -  \frac{i}{2} (\tau_3'-\tau_1') \, L_{2,\e} \notag  \\
&  \qquad
      -  \frac{1}{3}\,i\,(\tau_3''-\tau_1'')
      -  \frac{1}{3}\,i\,(\tau_3''-\tau_2'')
   \label{E-13}  \\
\wt{L}_{123,\e} - \wt{L}_{123,0}
&  =  - 2 \, \e^2 |\xi|^2 \sum_{j=1}^3 L_{j,\e} \, . \label{E-24}
\end{align}
\end{Lemma}

\begin{proof}
We have
\begin{align}
L_{1,\e} \circ L_{2,\e} \circ L_{3,\e}
&  =  (\partial_t - i\tau_{1,\e})
            \circ  (\partial_t - i\tau_{2,\e})
            \circ  (\partial_t -  i\tau_{3,\e}) \notag  \\
&  =  L_{123,\e}
      + (-i\tau_3') L_{1,\e} + \partial_t L_{23,\e}
      + (-i\tau_3'') \notag  \\
&  =  L_{123,\e}
      + (-i\tau_3') L_{1,\e}
      + (-i\tau_2') L_{3,\e}
      + (-i\tau_3') L_{2,\e} + (-i\tau_3'') \, . \label{E-12}
\end{align}
Summing over all permutations
we get~\eqref{E-tL123}.

Identity \eqref{E-13} follows from~\eqref{E-12}.

As \eqref{E-tL123} with $\e=0$ gives
\begin{equation} \label{E-tL123e0}
\wt{L}_{123,0}
  =  L
      +  \frac{1}{2} \, \partial_t \partial_\tau L
      +  \frac{1}{6} \, \partial_t^2 \partial_\tau^2 L \, ,
\end{equation}
we get
\begin{align*}
\wt{L}_{123,\e} - \wt{L}_{123,0}
&  =  L_{123,\e} - L
      +  \frac{1}{2} \, \partial_t \partial_\tau
         [ L_{123,\e} - L]
      +  \frac{1}{6} \, \partial_t^2 \partial_\tau^2
         [ L_{123,\e} - L ]  \\
&  =  - \e^2 |\xi|^2 \partial_\tau^2 L
   =  - \e^2 |\xi|^2 \partial_\tau^2 L_{123,\e}
   =  - 2 \, \e^2 |\xi|^2 \sum_{j=1}^3 L_{j,\e} \, . \qedhere
\end{align*}
\end{proof}

We define an energy,
after a Fourier transform with respect to the variables $x$ ($v= \mathcal {F}_x u$):
\[
E(t,\xi)
  \defeq  k(t,\xi)
      \Biggl[ \ \sum_{(j,h)\in\SS_2} |\wt{L}_{jh,\e}v|^2
              +  \HH^2(t,\xi) \,
                 \biggl[ \ \sum_{j=1}^3 |L_{j,\e} v|^2 + |v|^2 \ \biggr] \Biggr] ,
\]
where the weight $k(t,\xi)$ is defined by
\[
k(t,\xi)
  \defeq  \exp \Biggl[ -\eta \int_{-T}^t \KK(s,\xi) \, ds\Biggr] \, ,
\]
with
\begin{align*}
\KK(t,\xi)
  \defeq
&        \sum_{(j,h)\in\SS_2}
         \frac{|\tau_{j,\e}'-\tau_{h,\e}'|}
              {|\tau_{j,\e}-\tau_{h,\e}|}
         +  \sum_{(j,h)\in\SS_2}
            \frac{|\tau_{j,\e}''-\tau_{h,\e}''|}
                 {|\tau_{j,\e}' - \tau_{h,\e}'|+1}  \\
&  \quad  + \sum_{j=1}^3
            \frac{\bigl|\partial_t \widecheck{M}(\tau_{j,\e})\bigr|}
                 {\bigl|\widecheck{M}(\tau_{j,\e})\bigr|+1}
          + \sum_{j=1}^2
            \frac{\bigl|\partial_t \widecheck{N}(\sigma_{j,\e})\bigr|}
                 {\bigl|\widecheck{N}(\sigma_{j,\e})\bigr|+1}  \\
&  \quad  + \sum_{(j,h,l)\in\SS_3}
            \frac{\bigl|\widecheck{M}(\tau_{j,\e})\bigr|}
                 {|\tau_{j,\e}-\tau_{h,\e}|\cdot|\tau_{j,\e} - \tau_{l,\e}|}
          + \sum_{j=1}^2 \sqrt{ \frac{\bigl|\widecheck{N}(\sigma_{1,\e})\bigr|}
                                {|\sigma_{2,\e} - \sigma_{1,\e}|} \,}
          + \log|\xi| \, ,  \\
\HH(t,\xi)
&  \defeq  1 + \sum_{(j,h)\in\SS_2} \frac{|\tau_{j,\e}' - \tau_{h,\e}'|}
                                           {|\tau_{j,\e} - \tau_{h,\e}|}  \\
&  \quad  +  \sum_{(j,h,l)\in\SS_3}
             \frac{\bigl|\widecheck{M}(\tau_{j,\e})\bigr|}
                  {|\tau_{j,\e} - \tau_{h,\e}| \cdot
                   |\tau_{j,\e} - \tau_{l,\e}|}
          +  \sum_{j=1}^2 \sqrt{ \frac{\bigl|\widecheck{N}(\sigma_{j,\e})\bigr|+1}
                                 {|\sigma_{2,\e} - \sigma_{1,\e}|} \,} \, ,
\end{align*}
$\widecheck{M}$ and $\widecheck{N}$ are defined in~\eqref{E-cM} and~\eqref{E-cN}.

\smallskip

Differentiating the energy with respect to time we get
\begin{align*}
E'(t,\xi)
&  =  -\eta \, \KK(t,\xi) \, E(t,\xi)  \\
&  \qquad
      +  k(t,\xi)
      \Biggl[ \ 2 \, \sum_{(j,h)\in\SS_2}
                  \Re\< \partial_t \wt{L}_{jh,\e}v \, , \, \wt{L}_{jh,\e}v \>  \\
&  \qquad\qquad\qquad\qquad
              +  2 \, \HH(t,\xi) \,\HH'(t,\xi) \,
                 \biggl[ \ \sum_{j=1}^3 |L_{j,\e} v|^2 + |v|^2 \ \biggr]  \\
&  \qquad\qquad\qquad\qquad
              +  \HH^2(t,\xi) \,
                 \biggl[ \  \sum_{j=1}^3 2 \, \Re\<\partial_t L_{j,\e} v \, , \, L_{j,\e} v\>
                         + 2 \, \Re\<\partial_t v \, , \, v\> \ \biggr] \Biggr] \, .
\end{align*}

Now we show that the second, third and forth summand
can be estimated by $C\,\KK(t,\xi) \, E(t,\xi)$,
for some suitable positive constant $C$.

\subsection{Estimation of~the terms \
   $2\Re\< \partial_t \wt{L}_{jh,\e}v \, , \, \wt{L}_{jh,\e}v \>$}
\text{}

As
\[
\partial_t \wt{L}_{jh,\e}v
  =  (L_l \circ \wt{L}_{jh,\e}) v
      +  i\,\tau_{l,\e} \, \wt{L}_{jh,\e} v
\]
we have
\begin{align*}
2 \, \Re\< \partial_t \wt{L}_{jh,\e}v \, , \, \wt{L}_{jh,\e}v \>
&   =   2 \, \Re\< (L_l \circ \wt{L}_{jh,\e}) v \, , \, \wt{L}_{jh,\e}v \>
        +  2 \, \Re\< i\,\tau_{l,\e} \, \wt{L}_{jh,\e} v \, , \, \wt{L}_{jh,\e}v \>  \\
&   =   2 \, \Re\< (L_l \circ \wt{L}_{jh,\e}) v \, , \, \wt{L}_{jh,\e}v \> \, .
\end{align*}

Define
\begin{equation}
\begin{split}
\wt{M}
&  \defeq  \widecheck{M} + \frac{1}{2} \partial_t\partial_\tau \widecheck{M}  \\
&  =  M - \frac{1}{2} \partial_t\partial_\tau L
      + \frac{1}{2} \partial_t\partial_\tau M
      - \frac{1}{4} \partial_t^2\partial_\tau^2 L \, .
\end{split}
\end{equation}
so that
\begin{equation} \label{E-dec}
\wt{L}_{123,0} + \wt{M} + \widecheck{N}
  =  L + M + N \, ,
\end{equation}
hence
\begin{align*}
2 \, \Re\< (L_l \circ \wt{L}_{jh,\e}) v \, , \, \wt{L}_{jh,\e}v \>
&  =  2 \, \Re\< (L_l \circ \wt{L}_{jh,\e}) v - \wt{L}_{123,\e} v \, , \, \wt{L}_{jh,\e}v \>  \\
&  \qquad
     +  2 \, \Re\< \wt{L}_{123,\e} v - \wt{L}_{123,0} v \, , \, \wt{L}_{jh,\e}v \>  \\
&  \qquad
     +  2 \, \Re\< Lv + Mv + Nv \, , \, \wt{L}_{jh,\e}v \>  \\
&  \qquad
     -  2 \, \Re\< \wt{M} v \, , \, \wt{L}_{jh,\e}v \>
     -  2 \, \Re\< \widecheck{N} v \, , \, \wt{L}_{jh,\e}v \>
\end{align*}

\subsubsection{Estimation of~\
   $2 \, \Re\< (L_l \circ \wt{L}_{jh,\e}) v - \wt{L}_{123,\e} v \, , \, \wt{L}_{jh,\e}v \>$}
\text{}

First of~all, we have
\[
\Bigl|\Re\< (L_l \circ \wt{L}_{jh,\e}) v - \wt{L}_{123,\e} v \, , \, \wt{L}_{jh,\e}v \>\Bigr|
  \lesssim \bigl|(L_l \circ \wt{L}_{jh,\e}) v - \wt{L}_{123,\e} v\bigr| \,
            |\wt{L}_{jh,\e}v| \, .
\]

According to~\eqref{E-13},
$(L_l \circ \wt{L}_{jh,\e}) v - \wt{L}_{123,\e} v$
is a linear combination, with constant coefficients,
of terms like $(\tau_{\alpha,\e}'-\tau_{\beta,\e}') \, L_{\gamma,\e}v$,
with $(\alpha,\beta,\gamma)\in\SS_3$,
and $(\tau_{\alpha,\e}''-\tau_{\beta,\e}'')v$,
with~$(\alpha,\beta)\in\SS_2$,
hence:
\[
\bigl|(L_l \circ \wt{L}_{jh,\e}) v - \wt{L}_{123,\e} v\bigr|
  \lesssim  \sum_{(\alpha,\beta,\gamma)\in\SS_3}
            |\tau_{\alpha,\e}'-\tau_{\beta,\e}'| \, |L_{\gamma,\e}v|
            +  \sum_{(\alpha,\beta)\in\SS_2}
               |\tau_{\alpha,\e}''-\tau_{\beta,\e}''| \, |v| \, .
\]

\goodbreak

Concerning the first sum, from~\eqref{E-31} and~\eqref{E-32}, we have
\begin{align*}
L_{\gamma ,\e}v
&  =  i\,\frac{\wt{L}_{\gamma \alpha,\e}v-\wt{L}_{\gamma \beta,\e}v}
              {\tau_{\alpha,\e}-\tau_{\beta,\e}}
      -  \frac{1}{2} \,
         \frac{\tau_{\alpha,\e}' - \tau_{\beta,\e}'}
              {\tau_{\alpha,\e}-\tau_{\beta,\e}} \, v  \\
&  =  i\,\frac{\wt{L}_{\gamma \alpha,\e}v-\wt{L}_{\gamma \beta,\e}v}
              {\tau_{\alpha,\e}-\tau_{\beta,\e}}
      -  \frac{i}{2} \,
         \frac{\tau_{\alpha,\e}' - \tau_{\beta,\e}'}
              {\tau_{\alpha,\e}-\tau_{\beta,\e}} \,
         \frac{L_{\alpha,\e}v - L_{\beta,\e}v}
              {\tau_{\alpha,\e}-\tau_{\beta,\e}} \, ,
\intertext{%
hence}
|\tau_{\alpha,\e}'-\tau_{\beta,\e}'| \, |L_{\gamma,\e}v|
&  \lesssim  \biggl|\frac{\tau_{\alpha,\e}'-\tau_{\beta,\e}'}
                         {\tau_{\alpha,\e}-\tau_{\beta,\e}}\biggr|
             \bigl[|\wt{L}_{\gamma \alpha,\e}v| + |\wt{L}_{\gamma \beta,\e}v|\bigr]  \\
&  \qquad  +  \biggl|\frac{\tau_{\alpha,\e}'-\tau_{\beta,\e}'}
                          {\tau_{\alpha,\e}-\tau_{\beta,\e}}\biggr|^2 \,
              \bigl[|L_{\alpha,\e}v| + |L_{\beta,\e}v|\bigr]  \\
&  \lesssim  \KK(t,\xi) \,
             \biggl[ \sum_{(j,h)\in\SS_2} |\wt{L}_{jh,\e}v|
                     +  \HH(t,\xi) \, \sum_{j=1}^3 |L_{j,\e} v|\biggr] \, .
\end{align*}

Concerning the second sum,
from~\eqref{E-32},
we have
\begin{align*}
\bigl|(\tau_{\alpha,\e}'' - \tau_{\beta,\e}'' )v\bigr|
&  =  \frac{|\tau_{\alpha,\e}'' - \tau_{\beta,\e}''|}
           {|\tau_{\alpha,\e} - \tau_{\beta,\e}|} \,
      \bigl|L_{\beta,\e}v - L_{\alpha,\e}v\bigr|  \\
&  =  \frac{|\tau_{\alpha,\e}'' - \tau_{\beta,\e}''|}
           {|\tau_{\alpha,\e}' - \tau_{\beta,\e}'| + 1} \,
      \frac{|\tau_{\alpha,\e}' - \tau_{\beta,\e}'| + 1}
           {|\tau_{\alpha,\e} - \tau_{\beta,\e}|} \,
      \bigl|L_{\beta,\e}v - L_{\alpha,\e}v\bigr|  \\
&  \le  \KK(t,\xi) \, \HH(t,\xi) \, \sum_{j=1}^3 |L_{j,\e} v| \, .
\end{align*}

Combining the above estimates, we get
\[
\Bigl|\Re\< (L_l \circ \wt{L}_{jh,\e}) v - \wt{L}_{123,\e} v \, , \, \wt{L}_{jh,\e}v \>\Bigr|
  \lesssim  \KK(t,\xi) \,
             \biggl[ \sum_{(j,h)\in\SS_2} |\wt{L}_{jh,\e}v|^2
                     +  \HH^2(t,\xi) \, \sum_{j=1}^3 |L_{j,\e} v|^2\biggr] \, .
\]

\subsubsection{Estimation of~\
   $2 \, \Re\< \wt{L}_{123,\e} v - \wt{L}_{123,0} v \, , \, \wt{L}_{jh,\e}v \>$}
\text{}

Using~\eqref{E-24},
we have
\[
2 \, \Re\< \wt{L}_{123,\e} v - \wt{L}_{123,0} v \, , \, \wt{L}_{jh,\e}v \>
  \lesssim  \sum_{j=1}^3 |L_{j,\e}v|^2 + |\wt{L}_{jh,\e}v|^2 \, .
\]

\subsubsection{Estimation of~\
   $2 \, \Re\< Lv + Mv + Nv \, , \, \wt{L}_{jh,\e}v \>$}
\text{}

As
\[
Lv + Mv + Nv
  =  -pv
\]

\[
2 \, \Re\< Lv + Mv + Nv \, , \, \wt{L}_{jh,\e}v \>
  =  - 2 \, \Re\< pv \, , \, \wt{L}_{jh,\e}v \>
  \lesssim  |v|^2 + |\wt{L}_{jh,\e}v|^2 \, .
\]

\subsubsection{Estimation of~\
   $2 \, \Re\< \wt{M} v \, , \, \wt{L}_{jh,\e}v \>$}
\text{}

To estimate $\widecheck{M} = M - \frac{1}{2} \partial_t\partial_\tau L$
we can use the Lagrange interpolation formula
\begin{equation} \label{E-Lagrange}
\widecheck{M}
  =  \sum_{(j,h,l)\in\SS_3} \ell_j(t,\xi) \, L_{hl,\e} \, ,
\quad\text{where}\quad
\ell_j(t,\xi)
  \defeq  \frac{\widecheck{M}(\tau_{j,\e})}
          {(\tau_{j,\e}-\tau_{h,\e})(\tau_{j,\e}-\tau_{l,\e})} \, .
\end{equation}

If $\wt{M}$ were a linear combination of~the $\wt{L}_{hl,\e}$
then we could estimate~it.

Now we observe that from $\widecheck{M} = \sum \ell_j L_{hl,\e}$
it~does not follow $\wt{M} = \sum \ell_j \wt{L}_{hl,\e}$,
but it follows
\begin{align*}
\wt{M}
&  =  \widecheck{M} + \frac{1}{2} \partial_t\partial_\tau \widecheck{M}  \\
&  =  \sum_{(j,h,l)\in\SS_3} \ell_j L_{hl,\e}
      +  \frac{1}{2} \, \sum_{(j,h,l)\in\SS_3} \ell_j \partial_t\partial_\tau L_{hl,\e}
      +  \frac{1}{2} \, \sum_{(j,h,l)\in\SS_3} \partial_t \ell_j \partial_\tau L_{hl,\e}  \\
&  =  \sum_{(j,h,l)\in\SS_3} \ell_j \wt{L}_{hl,\e}
      +  \frac{1}{2} \,
         \sum_{(j,h,l)\in\SS_3} \partial_t \ell_j [L_{h,\e} + L_{l,\e}] \, .
\end{align*}

We have
\begin{align*}
\partial_t \ell_j
&  =  \partial_t
      \frac{\widecheck{M}(\tau_{j,\e})}
           {(\tau_{j,\e}-\tau_{h,\e}) \, (\tau_{j,\e}-\tau_{l,\e})}   \\
&  =  \frac{\partial_t\widecheck{M}(\tau_{j,\e})}
           {(\tau_{j,\e}-\tau_{h,\e}) \, (\tau_{j,\e}-\tau_{l,\e})}
      -  \frac{\widecheck{M}(\tau_{j,\e})}
              {(\tau_{j,\e}-\tau_{h,\e}) \, (\tau_{j,\e}-\tau_{l,\e})} \cdot
         \left( \frac{\tau_{j,\e}'-\tau_{h,\e}'}
                     {\tau_{j,\e}-\tau_{h,\e}}
                +  \frac{\tau_{j,\e}'-\tau_{l,\e}'}
                        {\tau_{j,\e}-\tau_{l,\e}}
              \right) ,
\end{align*}
and
\begin{align*}
\Biggl|\frac{\partial_t\widecheck{M}(\tau_{j,\e})}
           {(\tau_{j,\e}-\tau_{h,\e}) \, (\tau_{j,\e}-\tau_{l,\e})}\Biggr|
&  =  \frac{\bigl|\partial_t\widecheck{M}(\tau_{j,\e})\bigr|}
           {\bigl|\widecheck{M}(\tau_{j,\e})\bigr|+1} \cdot
      \frac{\bigl|\widecheck{M}(\tau_{j,\e})\bigr|+1}
           {|\tau_{j,\e}-\tau_{h,\e}|\,|\tau_{j,\e}-\tau_{l,\e}|}  \\
&  \lesssim  \frac{\bigl|\partial_t\widecheck{M}(\tau_{j,\e})\bigr|}
           {\bigl|\widecheck{M}(\tau_{j,\e})\bigr|+1}
        \left( \frac{\bigl|\widecheck{M}(\tau_{j,\e})\bigr|}
                    {|\tau_{j,\e}-\tau_{h,\e}|\,|\tau_{j,\e}-\tau_{l,\e}|}
                                 + 1 \right)
\end{align*}
since $|\tau_{j,\e}-\tau_{h,\e}|\ge C>0$.
Hence
\begin{equation} \label{E-deM}
\Biggl|\partial_t\frac{\widecheck{M}(\tau_{j,\e})}
           {(\tau_{j,\e}-\tau_{h,\e}) \, (\tau_{j,\e}-\tau_{l,\e})}\Biggr|
  \le  \KK(t,\xi) \, \HH(t,\xi) \, .
\end{equation}
Thus we get
\begin{align*}
|\wt{M} v|
&  \lesssim  \KK(t,\xi)
      \Biggl[ \ \sum_{(j,h)\in\SS_2} |\wt{L}_{jh,\e}v|
              +  \HH(t,\xi) \, \sum_{j=1}^3 |L_{j,\e} v| \Biggr] ,
\intertext{%
hence}
\biggl|2 \, \Re\< \wt{M} v \, , \, \wt{L}_{jh,\e}v \>\biggr|
&  \lesssim  \KK(t,\xi)
      \Biggl[ \ \sum_{(j,h)\in\SS_2} |\wt{L}_{jh,\e}v|^2
              +  \HH^2(t,\xi) \, \sum_{j=1}^3 |L_{j,\e} v|^2 \Biggr] .
\end{align*}

\subsubsection{Estimation of~\
   $2 \, \Re\< \widecheck{N} v \, , \, \wt{L}_{jh,\e}v \>$}
\text{}

Using Lagrange's interpolation formula:
\[
\widecheck{N}(\tau)
  =  \frac{\widecheck{N}(\sigma_{1,\e})}
          {\sigma_{1,\e} - \sigma_{2,\e}} \, (\tau-\sigma_{2,\e})
    +  \frac{\widecheck{N}(\sigma_{2,\e})}
            {\sigma_{2,\e} - \sigma_{1,\e}} \, (\tau-\sigma_{1,\e}) \, .
\]
Now, since $\tau_{1,\e} \le \sigma_{1,\e} \le \tau_{2,\e} \le \sigma_{2,\e} \le \tau_{3,\e}$,
we can find $\theta_1,\theta_2\in[0,1]$ such that
\[
\sigma_{j,\e}
  =  \theta_j \, \tau_{1,\e} + (1 - \theta_j) \, \tau_{3,\e} \, ,
\]
hence
\[
(\tau-\sigma_{j,\e})
  =  \theta_j \, (\tau-\tau_{1,\e}) + (1 - \theta_j) \, (\tau-\tau_{3,\e}) \, ,
\]
thus $\widecheck{N}(\tau)$ is a linear combination,
with bounded coefficients, of~terms of~the form
\[
\frac{\widecheck{N}(\sigma_{j,\e})}
     {\sigma_{1,\e} - \sigma_{2,\e}} \, (\tau-\tau_{h,\e}) \, ,
\]
$j=1,2$, $h=1,3$.

We get
\begin{align*}
2 \, \Re\< \widecheck{N} v \, , \, \wt{L}_{jh,\e}v \>
&  \lesssim
      \frac{\bigl|\widecheck{N}(\sigma_{1,\e})\bigr|
                    + \bigl|\widecheck{N}(\sigma_{2,\e})\bigr|}
                   {|\sigma_{2,\e} - \sigma_{1,\e}|}
    \Bigl[ |L_{1,\e}v| + |L_{3,\e}v| \Bigr] \, |\wt{L}_{jh,\e}v|  \\
&  \lesssim
    \sqrt{ \frac{\bigl|\widecheck{N}(\sigma_{1,\e})\bigr|
                 + \bigl|\widecheck{N}(\sigma_{2,\e})\bigr|}
                {|\sigma_{2,\e} - \sigma_{1,\e}|} \,}
    \Biggl[ |\wt{L}_{jh,\e}v|^2  \\
&  \qquad\qquad
       +  \frac{\bigl|\widecheck{N}(\sigma_{1,\e})\bigr|
                 + \bigl|\widecheck{N}(\sigma_{2,\e})\bigr|}
               {|\sigma_{2,\e} - \sigma_{1,\e}|}
          \Bigl[ |L_{1,\e}v|^2 + |L_{3,\e}v|^2 \Bigr] \Biggr]  \\
&  \lesssim  \KK(t,\xi)
      \Biggl[ \ \sum_{(j,h)\in\SS_2} |\wt{L}_{jh,\e}v|^2
              +  \HH^2(t,\xi) \, \sum_{j=1}^3 |L_{j,\e} v|^2 \Biggr] .
\end{align*}

\subsection{Estimation of~the term \
   $2 \, \HH(t,\xi) \,\HH'(t,\xi) $}
\text{}

\begin{Lemma}
There exists $C_H>0$ such that
\[
\HH'(t,\xi)
  \le  C_H \, \KK(t,\xi) \, \HH(t,\xi) \, .
\]
\end{Lemma}

\begin{proof}
For any $(j,h)\in\SS_2$ we have
\begin{align*}
\partial_t \frac{\tau_{j,\e}' - \tau_{h,\e}'}
                {\tau_{j,\e} - \tau_{h,\e}}
&  =  \frac{\tau_{j,\e}'' - \tau_{h,\e}''}
           {\tau_{j,\e} - \tau_{h,\e}}
      -
      \frac{(\tau_{j,\e}' - \tau_{h,\e}')^2}
           {(\tau_{j,\e} - \tau_{h,\e})^2} \, ,
\intertext{%
hence}
\partial_t \frac{|\tau_{j,\e}' - \tau_{h,\e}'|}
                {|\tau_{j,\e} - \tau_{h,\e}|}
&  \le  \frac{|\tau_{j,\e}'' - \tau_{h,\e}''|}
             {|\tau_{j,\e}' - \tau_{h,\e}'|+1}
      \cdot
      \frac{|\tau_{j,\e}' - \tau_{h,\e}'|+1}
           {|\tau_{j,\e} - \tau_{h,\e}|}
      +
      \frac{|\tau_{j,\e}' - \tau_{h,\e}'|^2}
           {|\tau_{j,\e} - \tau_{h,\e}|^2}  \\
&  \lesssim  \KK(t,\xi) \, \HH(t,\xi) \, .
\end{align*}

The terms $\partial_t\frac{\widecheck{M}(\tau_{j,\e})}
           {(\tau_{j,\e}-\tau_{h,\e}) \, (\tau_{j,\e}-\tau_{l,\e})}$
are estimated as in~\eqref{E-deM}.

For any $(j,h)\in\SS_2$ we have
\begin{align*}
\partial_t \frac{\bigl|\widecheck{N}(\sigma_{j,\e})\bigr|+1}
                {|\sigma_{2,\e} - \sigma_{1,\e}|}
&  =  \frac{1}{|\sigma_{2,\e} - \sigma_{1,\e}|}
      \cdot
      \frac{\widecheck{N}(\sigma_{j,\e})}
           {\bigl|\widecheck{N}(\sigma_{j,\e})\bigr|}
      \cdot
      \partial_t \widecheck{N}(\sigma_{j,\e})
   -  \frac{\bigl|\widecheck{N}(\sigma_{j,\e})\bigr|+1}
           {|\sigma_{2,\e} - \sigma_{1,\e}|}
      \cdot
      \frac{\sigma_{2,\e}' - \sigma_{1,\e}'}
           {\sigma_{2,\e} - \sigma_{1,\e}}  \\
&  =  \frac{\bigl|\widecheck{N}(\sigma_{j,\e})\bigr| + 1}
          {|\sigma_{2,\e} - \sigma_{1,\e}|}
      \Biggl[
      \frac{\widecheck{N}(\sigma_{j,\e})}
           {\bigl|\widecheck{N}(\sigma_{j,\e})\bigr|}
      \cdot
      \frac{\partial_t \widecheck{N}(\sigma_{j,\e})}
           {\bigl|\widecheck{N}(\sigma_{j,\e})\bigr| + 1}
   -  \frac{\sigma_{2,\e}' - \sigma_{1,\e}'}
           {\sigma_{2,\e} - \sigma_{1,\e}}
      \Biggr] \, ,
\end{align*}
hence
\begin{align*}
\partial_t \sqrt{\frac{\bigl|\widecheck{N}(\sigma_{j,\e})\bigr|+1}
                {|\sigma_{2,\e} - \sigma_{1,\e}|}\,}
&  =  \frac{1}{2\sqrt{\dfrac{\bigl|\widecheck{N}(\sigma_{j,\e})\bigr|+1}
                            {|\sigma_{2,\e} - \sigma_{1,\e}|}}\,} \,
      \partial_t \frac{\bigl|\widecheck{N}(\sigma_{j,\e})\bigr| + 1}
                      {|\sigma_{2,\e} - \sigma_{1,\e}|}  \\
&  =  \frac{1}{2} \,
      \sqrt{\frac{\bigl|\widecheck{N}(\sigma_{j,\e})\bigr|+1}
                 {|\sigma_{2,\e} - \sigma_{1,\e}|}\,}
      \Biggl[
      \frac{\widecheck{N}(\sigma_{j,\e})}
           {\bigl|\widecheck{N}(\sigma_{j,\e})\bigr|}
      \cdot
      \frac{\partial_t \widecheck{N}(\sigma_{j,\e})}
           {\bigl|\widecheck{N}(\sigma_{j,\e})\bigr| + 1}
   -  \frac{\sigma_{2,\e}' - \sigma_{1,\e}'}
           {\sigma_{2,\e} - \sigma_{1,\e}}
      \Biggr] .
\end{align*}

All the terms but the last can be estimated by $\KK(t,\xi)$ or $\HH(t,\xi)$.
For the last, we remark that (cf.~Lemma~\ref{L-A-1})
\[
(\tau_{1,\e}-\tau_{2,\e})^2+(\tau_{2,\e}-\tau_{3,\e})^2+(\tau_{3,\e}-\tau_{1,\e})^2
  =  \frac{9}{2} \, (\sigma_{2,\e}-\sigma_{1,\e})^2 \, ,
\]
\goodbreak
thus
\begin{align*}
\frac{\sigma_{2,\e}' - \sigma_{1,\e}'}{\sigma_{2,\e} - \sigma_{1,\e}}
&  =  \frac{1}{2} \,
      \frac{[(\sigma_{2,\e} - \sigma_{1,\e})^2]'}{(\sigma_{2,\e} - \sigma_{1,\e})^2}
   =  \frac{1}{2} \,
      \frac{\Bigl[ \sum_{j,h\in\SS_2} (\tau_{j,\e}-\tau_{h,\e})^2 \Bigr]'}
                   {\sum_{j,h\in\SS_2} (\tau_{j,\e}-\tau_{h,\e})^2}  \\
&  =  \sum_{j,h\in\SS_2}
      \frac{(\tau_{j,\e}-\tau_{h,\e}) \, (\tau_{j,\e}'-\tau_{h,\e}') }
           {\sum_{j,h\in\SS_2} (\tau_{j,\e}-\tau_{h,\e})^2} \, ,
\intertext{%
which gives}
\biggl|\frac{\sigma_{2,\e}' - \sigma_{1,\e}'}{\sigma_{2,\e} - \sigma_{1,\e}} \biggr|
&  \le  \sum_{j,h\in\SS_2} \frac{|\tau_{j,\e}'-\tau_{h,\e}'|}{|\tau_{j,\e}-\tau_{h,\e}|} \, ,
\end{align*}
which can be estimated by $\KK(t,\xi)$ or $\HH(t,\xi)$.

Finally we get
\[
\partial_t \sqrt{\frac{\bigl|\widecheck{N}(\sigma_{j,\e})\bigr|+1}
                {|\sigma_{2,\e} - \sigma_{1,\e}|}\,}
  \lesssim  \KK(t,\xi) \, \HH(t,\xi) \, . \qedhere
\]
\end{proof}

\subsection{Estimation of~the terms \
   $2 \, \HH^2 \Re\< \partial_t L_{j,\e}v \, , \, L_{j,\e}v \>$}
\text{}

As
\begin{align*}
\partial_t (L_{j,\e} v)
&  =  [L_{h,\e} + i \tau_{h,\e}] L_{j,\e} v  \\
&  =  (L_{h,\e} \circ L_{j,\e}) v + i \tau_{h,\e} L_{j,\e} v \, ,
\end{align*}
we have
\begin{align*}
2\Re\< \partial_t L_{j,\e}v \, , \, L_{j,\e}v \>
&  =  2\Re\< (L_{h,\e} \circ L_{j,\e}) v \, , \, L_{j,\e}v \>
      +  2\Re\< i \tau_{h,\e} L_{j,\e} \, , \, L_{j,\e}v \>  \\
&  =  2\Re\< (L_{h,\e} \circ L_{j,\e}) v \, , \, L_{j,\e}v \> \, ,
\intertext{%
as $\tau_{h,\e}$ is a real function,
hence}
2 \, \HH^2 \Re\<\partial_t L_{j,\e} v \, , \, L_{j,\e} v\>
&  \le  \KK \, \Bigl[ \bigl|(L_{h,\e} \circ L_{j,\e}) v\bigr|^2
                      +  \HH^2 \, |L_{j,\e}|^2 \Bigr] \, .
\end{align*}

Now, using~\eqref{E-23}, we have
\begin{align*}
L_{j,\e} \circ L_{h,\e} - \wt{L}_{jh,\e}
&  =  \frac{i}{2} \, (\tau_{j,\e}' - \tau_{h,\e}')  \\
&  =  -\frac{i}{2} \,
       \frac{\tau_{j,\e}' - \tau_{h,\e}'}
            {\tau_{j,\e} - \tau_{h,\e}} (L_{j,\e}v - L_{h,\e}v) \, ,
\end{align*}
hence
\begin{align*}
\bigl|(L_{h,\e} \circ L_{j,\e}) v\bigr|^2
&  \le  2\bigl|\wt{L}_{jh,\e}v\bigr|^2
        +2\bigl|(L_{j,\e} \circ L_{h,\e} - \wt{L}_{jh,\e}) v\bigr|^2  \\
&  \le  2\bigl|\wt{L}_{jh,\e}v\bigr|^2
        +  \Bigl|\frac{\tau_{j,\e}' - \tau_{h,\e}'}
                      {\tau_{j,\e} - \tau_{h,\e}}\Bigr|^2
           \Bigl(|L_{j,\e}v|^2 + |L_{h,\e}v|^2\Bigr) \, .
\end{align*}

\subsection{Estimation of~the terms \
   $2\Re\< \partial_t v \, , \, v \>$}
\text{}

As
\[
\partial_t v
  =  L_{1,\e} v + i \tau_{1,\e} v \, ,
\]
we have
\begin{align*}
2 \, \Re\<\partial_t v \, , \, v\>
&  =  2 \, \Re\<L_{1,\e} v \, , \, v\>
      +  2 \Re\<i \tau_{1,\e} v \, , \, v\>  \\
&  \le  |L_{1,\e} v|^2 + |v|^2 \, .
\end{align*}

Taking $\eta$ large we arrive to the estimates
\[
E'(t)
  \le  C \, E(t) \, ,
\]

\[
E(t)
  \le  \exp\bigl( C(t-t_0)\bigr) \, E(t_0) \, ,
\]
{}from this taking into account the inequality
\[
\KK(t,\xi) \ge \bigl(2+|\xi|\bigr)^{-c_0}
\]
it~follows that the Cauchy problem for the given equation is well posed.

This concludes the proof of~Theorem~\ref{T-2}.

\section{Equivalent forms of~the Levi conditions} \label{S-EF}

In this paragraph we can give some alternative forms of~the Levi conditions.
In particular we express these conditions in terms
of the roots of~$L$.

We recall that
\begin{align}
\bigl|\tau_{j,\e}-\tau_{h,\e}\bigl|
&  \approx  \bigl|\tau_j-\tau_h\bigl|+1 \, ,  &
&  \text{for any $(j,h)\in\SS_2$} \, , \label{E-27}
\intertext{%
and}
\bigl|\tau_j-\tau_{j,\e}\bigl|
&  \lesssim  1 \, ,  &
&  \text{for any $j=1,2,3$} \, , \label{E-28}
\end{align}

\begin{Proposition} \label{P-M}
Hypothesis~\eqref{E-Mcl} is equivalent
to the conditions
\begin{subequations}
\begin{gather}
\int_0^T \sum_{(j,h,l)\in\SS_3} \frac{\bigl|\widecheck{M}\bigl(t,\tau_j(t,\xi),\xi\bigr)\bigr|}
              {\Bigl(\bigl|\tau_j(t,\xi)-\tau_h(t,\xi)\bigr|+1\Bigr)
               \Bigl(\bigl|\tau_j(t,\xi)-\tau_l(t,\xi)\bigr|+1\Bigr)} \, dt
  \lesssim \log\bigl(1+|\xi|\bigr) \, , \label{E-Mct}  \\
\int_0^T \frac{\bigl|\partial_\tau \widecheck{M}\bigl(t,\tau_1(t,\xi),\xi\bigr)\bigr|
               + \bigl|\partial_\tau \widecheck{M}\bigl(t,\tau_3(t,\xi),\xi\bigr)\bigr|}
              {\bigl|\tau_1(t,\xi)-\tau_3(t,\xi)\bigr|+1} \, dt
  \lesssim \log\bigl(1+|\xi|\bigr) \, . \label{E-dMct}
\end{gather}
\end{subequations}
\end{Proposition}

\begin{proof}
For the sake of~simplicity,
we omit the $t$ and $\xi$ variables.

\smallskip

We start by proving that~\eqref{E-Mcl} implies~\eqref{E-dMct}.

By the Lagrange interpolation formula we have
\[
\widecheck{M}(\tau)
  =  \sum_{(j,h,l)\in\SS_3} \ell_{j,\e} \, L_{hl,\e} (\tau) \, ,
\quad\text{where}\quad
\ell_{j,\e}
  \defeq  \frac{\widecheck{M}(\tau_{j,\e})}
          {(\tau_{j,\e}-\tau_{h,\e})(\tau_{j,\e}-\tau_{l,\e})} \, .
\]
Differentiating with respect to~$\tau$:
\begin{align*}
\partial_\tau \widecheck{M}(\tau)
&  =  \sum_{(j,h,l)\in\SS_3} \ell_{j,\e} \,
      \bigl[L_{k,\e}(\tau) + L_{l,\e}(\tau)\bigr]  \\
&  =  \sum_{(j,h,l)\in\SS_3} \bigl[ \ell_{j,\e} + \ell_{k,\e} \bigr] \, L_{l,\e}(\tau) \, ,
\end{align*}
hence $\partial_\tau \widecheck{M}(\tau)$ is a linear combination
of~$L_{1,\e}(\tau)$, $L_{2,\e}(\tau)$ and~$L_{3,\e}(\tau)$
with coefficients verifying the logaritmic condition.

To prove that $\partial_\tau \widecheck{M}(\tau)$ is a linear combination
of only $L_{1,\e}$ and~$L_{3,\e}$
we~note that,
since $\tau_{1,\e} \le \tau_{2,\e} \le \tau_{3,\e}$,
we can find $\theta\in[0,1]$ such that
\[
\tau_{2,\e}
  =  \theta \, \tau_{1,\e} + (1 - \theta) \, \tau_{3,\e} \, ,
\]
hence
\[
L_{2,\e}
  =  \theta \, L_{1,\e} + (1 - \theta) \, L_{3,\e} \, ,
\]
and then
\begin{equation} \label{E-62}
\partial_\tau \widecheck{M}(\tau)
  =  b_{1,\e} \, L_{1,\e}(\tau) + b_{3,\e} \, L_{3,\e}(\tau) \, ,
\end{equation}
where $b_{1,\e}$ and $b_{3,\e}$
are some linear combination of~the $\ell_{1,\e}$, $\ell_{2,\e}$ and~$\ell_{3,\e}$,
hence they verify the logaritmic condition.

Substituting $\tau$ with $\tau_{1,\e}$ and~$\tau_{3,\e}$
in~\eqref{E-62} we~get
\begin{equation} \label{E-60}
\int_0^T \frac{\bigl|\partial_\tau \widecheck{M}(\tau_{1,\e})\bigr|
                +  \bigl|\partial_\tau \widecheck{M}(\tau_{3,\e})\bigr|}
              {\bigl|\tau_{1,\e}-\tau_{3,\e}\bigl|} \, dt
  \lesssim \log\bigl(1+|\xi|\bigr) \, .
\end{equation}

Since
\begin{align*}
\partial_\tau \widecheck{M}(\tau_1)
&  =  \partial_\tau \widecheck{M}(\tau_{1,\e})
      +  \partial_\tau^2 \widecheck{M}(\tau_{1,\e}) \, (\tau_1 - \tau_{1,\e})
\intertext{%
we get}
\bigl|\partial_\tau \widecheck{M}(\tau_1)\bigr|
&  \lesssim  \bigl|\partial_\tau \widecheck{M}(\tau_{1,\e})\bigr| + 1 \, ,
\intertext{%
hence}
\frac{\bigl|\partial_\tau \widecheck{M}(\tau_1)\bigr|}
     {|\tau_3-\tau_1|+1}
&  \lesssim  \frac{\bigl|\partial_\tau \widecheck{M}(\tau_{1,\e})\bigr|}
                  {|\tau_{3,\e}-\tau_{1,\e}|} + 1 \, .
\end{align*}

An analogous estimate holds true with $\tau_1$ and $\tau_{1,\e}$
replaced by~$\tau_3$ and~$\tau_{3,\e}$.
Combining such inequalities with~\eqref{E-60} we get~\eqref{E-dMct}.

\smallskip

To prove that~\eqref{E-Mcl} implies~\eqref{E-Mct}
we~remark that, since $\widecheck{M}$ is a polynomial of~degree~$2$,
it~coincides with its Taylor's expansion of~order~$2$,
hence we have:
\[
\widecheck{M}(\tau_1)
  =  \widecheck{M}(\tau_{1,\e})
     +  \partial_\tau \widecheck{M}(\tau_{1,\e}) \, (\tau_1 - \tau_{1,\e})
     +  \frac{1}{2} \, \partial_\tau^2 \widecheck{M}(\tau_{1,\e}) \, (\tau_1 - \tau_{1,\e})^2 \, .
\]
Taking into account~\eqref{E-27} and~\eqref{E-28},
we have
\[
\frac{\bigl|\widecheck{M}(\tau_1)\bigr|}
     {\bigl(|\tau_2-\tau_1|+1\bigr) \, \bigl(|\tau_3-\tau_1|+1\bigr)}
  \lesssim  \frac{\bigl|\widecheck{M}(\tau_{1,\e})\bigr|}
                 {|\tau_{2,\e}-\tau_{1,\e}| \, |\tau_{3,\e}-\tau_{1,\e}|}
         +  \frac{\bigl|\partial_\tau \widecheck{M}(\tau_{1,\e})\bigr|}
                 {|\tau_{3,\e}-\tau_{1,\e}|} + 1 \, .
\]
An analogous estimate holds true for $\widecheck{M}(\tau_3)$.

To prove the estimate for $\widecheck{M}(\tau_2)$,
we split the phase space $[0,T] \times \R^n$
in two sub-zones:
\begin{align*}
Z_1
&  =  \Bigl\{ \ (t,\xi) \in [0,T] \times \R^n \ \Bigm|
              \ |\tau_2-\tau_1| \le |\tau_3-\tau_2| \ \Bigr\}  \\
Z_2
&  =  \Bigl\{ \ (t,\xi) \in [0,T] \times \R^n \ \Bigm|
              \ |\tau_3-\tau_2| \le |\tau_2-\tau_1| \ \Bigr\} \, .
\end{align*}
For $(t,\xi) \in Z_1$ we write
\[
\widecheck{M}(\tau_2)
  =  \widecheck{M}(\tau_{1,\e})
     +  \partial_\tau \widecheck{M}(\tau_{1,\e}) \, (\tau_2 - \tau_{1,\e})
     +  \frac{1}{2} \, \partial_\tau^2 \widecheck{M}(\tau_{1,\e}) \, (\tau_2 - \tau_{1,\e})^2 \, ,
\]
and, since
\[
|\tau_2 - \tau_{1,\e}|
  \lesssim  |\tau_2 - \tau_1| + 1
  \le  |\tau_3 - \tau_2| + 1 \, ,
\]
we get
\begin{align*}
&  \frac{\bigl|\widecheck{M}(\tau_2)\bigr|}
        {\bigl(|\tau_2-\tau_1|+1\bigr) \, \bigl(|\tau_3-\tau_2|+1\bigr)}  \\
&  \qquad
   \lesssim  \frac{\bigl|\widecheck{M}(\tau_{1,\e})\bigr|}
                  {\bigl(|\tau_2-\tau_1|+1\bigr) \, \bigl(|\tau_3-\tau_2|+1\bigr)}
           +  \frac{\bigl|\partial_\tau \widecheck{M}(\tau_{1,\e})\bigr|\,
                   |\tau_2 - \tau_{1,\e}|}
                  {\bigl(|\tau_2-\tau_1|+1\bigr) \, \bigl(|\tau_3-\tau_2|+1\bigr)}  \\
&  \qquad\qquad
           +  \frac{1}{2} \, \bigl|\partial_\tau^2 \widecheck{M}(\tau_{1,\e})\bigr| \,
              \frac{|\tau_2 - \tau_{1,\e}|^2}
                  {\bigl(|\tau_2-\tau_1|+1\bigr) \, \bigl(|\tau_3-\tau_2|+1\bigr)}  \\
&  \qquad
   \lesssim  \frac{\bigl|\widecheck{M}(\tau_{1,\e})\bigr|}
                  {\bigl(|\tau_2-\tau_1|+1\bigr) \, \bigl(|\tau_3-\tau_2|+1\bigr)}
           +  \frac{\bigl|\partial_\tau \widecheck{M}(\tau_{1,\e})\bigr|}
                   {|\tau_3-\tau_2|+1} + 1 \, .
\end{align*}

Now, using the fact that
$(t,\xi) \in Z_1$
if, and only if
$|\tau_3-\tau_1| \le 2 |\tau_3-\tau_2|$,
we~can estimate the second term by
$\frac{\bigl|\partial_\tau \widecheck{M}(\tau_{1,\e})\bigr|}
      {|\tau_3-\tau_1|+1}$.
Finally, by~\eqref{E-27},
we get
\begin{equation} \label{E-30}
\frac{\bigl|\widecheck{M}(\tau_2)\bigr|}
        {\bigl(|\tau_2-\tau_1|+1\bigr) \, \bigl(|\tau_3-\tau_2|+1\bigr)}
   \lesssim  \frac{\bigl|\widecheck{M}(\tau_{1,\e})\bigr|}
                  {|\tau_{2,\e}-\tau_{1,\e}| \, |\tau_{3,\e}-\tau_{1,\e}|}
          +  \frac{\bigl|\partial_\tau \widecheck{M}(\tau_{1,\e})\bigr|}
                  {|\tau_{3,\e}-\tau_{1,\e}|} + 1 \, .
\end{equation}

For $(t,\xi) \in Z_2$ we repeat the above calculation,
with $\tau_{1,\e}$ and $\tau_{3,\e}$ exhanged,
and we get:
\begin{equation} \label{E-33}
\frac{\bigl|\widecheck{M}(\tau_2)\bigr|}
        {\bigl(|\tau_2-\tau_1|+1\bigr) \, \bigl(|\tau_3-\tau_2|+1\bigr)}  \\
   \lesssim  \frac{\bigl|\widecheck{M}(\tau_{3,\e})\bigr|}
                  {|\tau_{3,\e}-\tau_{2,\e}| \, |\tau_{3,\e}-\tau_{1,\e}|}
          +  \frac{\bigl|\partial_\tau \widecheck{M}(\tau_{3,\e})\bigr|}
                  {|\tau_{3,\e}-\tau_{1,\e}|} + 1 \, .
\end{equation}

Combining~\eqref{E-30} and~\eqref{E-33} we get
\begin{multline*}
\frac{\bigl|\widecheck{M}(\tau_2)\bigr|}
        {\bigl(|\tau_2-\tau_1|+1\bigr) \, \bigl(|\tau_3-\tau_2|+1\bigr)}  \\
   \lesssim  \frac{\bigl|\widecheck{M}(\tau_{3,\e})\bigr|}
                  {|\tau_{3,\e}-\tau_{2,\e}| \, |\tau_{3,\e}-\tau_{1,\e}|}
          +  \frac{\bigl|\widecheck{M}(\tau_{1,\e})\bigr|}
                  {|\tau_{2,\e}-\tau_{1,\e}| \, |\tau_{3,\e}-\tau_{1,\e}|}
          +  \frac{\bigl|\partial_\tau \widecheck{M}(\tau_{1,\e})\bigr|
                   +  \bigl|\partial_\tau \widecheck{M}(\tau_{3,\e})\bigr|}
                  {|\tau_{3,\e}-\tau_{1,\e}|} + 1 \, .
\end{multline*}

Thanks to~\eqref{E-60} we see that~\eqref{E-Mcl} implies~\eqref{E-Mct}.

\smallskip

Finally we prove that~\eqref{E-Mct} and~\eqref{E-dMct} imply~\eqref{E-Mcl}.
Using Taylor expansion as before, we~have
\begin{equation} \label{E-Tay}
\widecheck{M}(\tau_{j,\e})
  =  \widecheck{M}(\tau_j)
     +  \partial_\tau \widecheck{M}(\tau_j) \, (\tau_{j,\e} - \tau_j)
     +  \frac{1}{2} \, \partial_\tau^2 \widecheck{M}(\tau_j) \, (\tau_{j,\e} - \tau_j)^2 \, ,
\quad
\text{for $j=1,2,3$} \, .
\end{equation}
Using~\eqref{E-Tay}
with $j=1$ we have
\begin{align*}
\frac{\bigl|\widecheck{M}(\tau_{1,\e})\bigr|}
     {|\tau_{2,\e}-\tau_{1,\e}|\,|\tau_{3,\e}-\tau_{1,\e}|}
&
  \lesssim  \frac{\bigl|\widecheck{M}(\tau_1)\bigr|}
                 {|\tau_{2,\e}-\tau_{1,\e}|\,|\tau_{3,\e}-\tau_{1,\e}|}
            +  \frac{\bigl|\partial_\tau \widecheck{M}(\tau_1)\bigr|\,|\tau_{1,\e} - \tau_1|}
                    {|\tau_{2,\e}-\tau_{1,\e}|\,|\tau_{3,\e}-\tau_{1,\e}|}  \\
&
  \qquad   +  \frac{|\tau_{1,\e} - \tau_1|^2}{|\tau_{2,\e}-\tau_{1,\e}|\,|\tau_{3,\e}-\tau_{1,\e}|} \, ,
\end{align*}
and, using~\eqref{E-27} and~\eqref{E-28}, we~get
\[
\frac{\bigl|\widecheck{M}(\tau_{1,\e})\bigr|}
     {|\tau_{2,\e}-\tau_{1,\e}|\,|\tau_{3,\e}-\tau_{1,\e}|}
  \lesssim  \frac{\bigl|\widecheck{M}(\tau_1)\bigr|}
                 {\bigl(|\tau_2-\tau_1|+1\bigl) \,\bigl(|\tau_3-\tau_1|+1\bigl)}
            +  \frac{\bigl|\partial_\tau \widecheck{M}(\tau_1)\bigr|}
                    {|\tau_3-\tau_1|+1}
            +  1 \, .
\]
Analogously, exchanging $\tau_1$ with $\tau_3$
and~$\tau_{1,\e}$ with $\tau_{3,\e}$:
\[
\frac{\bigl|\widecheck{M}(\tau_{3,\e})\bigr|}
     {|\tau_{2,\e}-\tau_{3,\e}|\,|\tau_{1,\e}-\tau_{3,\e}|}
  \lesssim  \frac{\bigl|\widecheck{M}(\tau_3)\bigr|}
                 {\bigl(|\tau_3-\tau_2|+1\bigl) \,\bigl(|\tau_3-\tau_1|+1\bigl)}
            +  \frac{\bigl|\partial_\tau \widecheck{M}(\tau_3)\bigr|}
                    {|\tau_3-\tau_1|+1}
            +  1 \, .
\]

Using~\eqref{E-Tay} with $j=2$, \eqref{E-27} and~\eqref{E-28}, we have:
\begin{align*}
\frac{\bigl|\widecheck{M}(\tau_{2,\e})\bigr|}
     {|\tau_{2,\e}-\tau_{1,\e}|\,|\tau_{3,\e}-\tau_{2,\e}|}
&  \lesssim  \frac{\bigl|\widecheck{M}(\tau_2)\bigr|}
                 {|\tau_{2,\e}-\tau_{1,\e}|\,|\tau_{3,\e}-\tau_{2,\e}|}
            +  \frac{\bigl|\partial_\tau \widecheck{M}(\tau_2)\bigr|\,|\tau_{2,\e} - \tau_2|}
                    {|\tau_{2,\e}-\tau_{1,\e}|\,|\tau_{3,\e}-\tau_{2,\e}|}  \\
&
  \qquad   +  \frac{|\tau_{2,\e} - \tau_2|^2}{|\tau_{2,\e}-\tau_{1,\e}|\,|\tau_{3,\e}-\tau_{2,\e}|}  \\
&
  \lesssim  \frac{\bigl|\widecheck{M}(\tau_2)\bigr|}
                 {\bigl(|\tau_3-\tau_2|+1\bigl) \,\bigl(|\tau_2-\tau_1|+1\bigl)}
            +  \frac{\bigl|\partial_\tau \widecheck{M}(\tau_2)\bigr|}
                    {\bigl(|\tau_3-\tau_2|+1\bigl) \,\bigl(|\tau_2-\tau_1|+1\bigl)}
            +  1 \, .
\end{align*}

To estimate the second term we~note that,
since $\tau_1 \le \tau_2 \le \tau_3$,
we can find $\theta\in[0,1]$ such that
\[
\tau_2
  =  \theta \, \tau_1 + (1 - \theta) \, \tau_3 \, ,
\]
hence
\[
\partial_\tau \widecheck{M}(\tau_2)
  =  \theta \, \partial_\tau \widecheck{M}(\tau_1)
     +  (1 - \theta) \, \partial_\tau \widecheck{M}(\tau_3) \, ,
\]
and then
\[
\bigl|\partial_\tau \widecheck{M}(\tau_2)\bigr|
  \le  \bigl|\partial_\tau \widecheck{M}(\tau_1)\bigr|
     +  \bigl|\partial_\tau \widecheck{M}(\tau_3)\bigr| \, .
\]
We note also that if $(t,\xi) \in Z_1$
then
$|\tau_3-\tau_1| \le 2 |\tau_3-\tau_2|$,
whereas if $(t,\xi) \in Z_2$
then
$|\tau_3-\tau_1| \le 2 |\tau_2-\tau_1|$,
thus
\[
\bigl(|\tau_3-\tau_2|+1\bigl) \,\bigl(|\tau_2-\tau_1|+1\bigl)
   \gtrsim  |\tau_3-\tau_1|+1 \, .
\]

Combining the above estimates, we get
\[
\frac{\bigl|\partial_\tau \widecheck{M}(\tau_2)\bigr|}
     {\bigl(|\tau_3-\tau_2|+1\bigl) \,\bigl(|\tau_2-\tau_1|+1\bigl)}
  \le  \frac{\bigl|\partial_\tau \widecheck{M}(\tau_1)\bigr|+\bigl|\partial_\tau \widecheck{M}(\tau_3)\bigr|}
            {|\tau_3-\tau_1|+1} \, . \qedhere
\]
\end{proof}

\begin{Proposition} \label{P-Ncpp}
Hypothesis~\eqref{E-Ncl} is equivalent to the condition
\begin{equation} \label{E-Ncpp}
\int_0^T \sum_{j=1}^2
  \sqrt{\frac{\bigl|\widecheck{N}\bigl(t,\sigma_j(t,\xi),\xi\bigr)\bigr|}
       {\bigl|\sigma_2(t,\xi) - \sigma_1(t,\xi)\bigr|+1}\,} \, dt
  \logestimate \, ,
\end{equation}
where $\sigma_1(t,\xi)$ and~$\sigma_2(t,\xi)$ are the roots in~$\tau$
of $\partial_\tau L(t,\tau,\xi)$.
\end{Proposition}

\begin{proof}
First of~all, we remark that, since
\[
\Bigl| \bigl|\widecheck{N}(\mu_1)\bigr| - \bigl|\widecheck{N}(\mu_2)\bigr|\Bigr|
  \lesssim  |\mu_2-\mu_1| \, ,
\]
condition~\eqref{E-Ncl} is equivalent to
\begin{equation} \label{E-Nc1p}
\int_0^T
     \sqrt{ \frac{\bigl|\widecheck{N}\bigl(t,\mu_{1}(t,\xi),\xi\bigr)\bigr|}
                 {\bigl|\mu_{2}(t,\xi) - \mu_{1}(t,\xi)\bigr|} \, } \, dt
  \logestimate \, .
\end{equation}

Analogously, since
\[
\Bigl| \bigl|\widecheck{N}(\sigma_1)\bigr| - \bigl|\widecheck{N}(\sigma_2)\bigr|\Bigr|
  \lesssim  |\sigma_2-\sigma_1| \, ,
\]
condition~\eqref{E-Ncpp} is equivalent to
\begin{equation} \label{E-Ncpp1}
\int_0^T
  \sqrt{\frac{\bigl|\widecheck{N}\bigl(t,\sigma_1(t,\xi),\xi\bigr)\bigr|}
       {\bigl|\sigma_2(t,\xi) - \sigma_1(t,\xi)\bigr|+1}\,} \, dt
  \logestimate \, .
\end{equation}

Now, by~\eqref{E-DDD}, we see that
\[
|\mu_2-\mu_1|
  \approx  |\sigma_2-\sigma_1|+1 \, ,
\]
whereas, by direct computation, we see that
\[
|\mu_1-\sigma_1|
  \lesssim 1 \, ,
\]
from which we deduce the equivalence of~\eqref{E-Nc1p} and~\eqref{E-Ncpp1}.
\end{proof}

\begin{Remark} \label{R-33}
More generally,
Hypothesis~\eqref{E-Ncl} is equivalent to each of~the conditions
\begin{equation} \label{E-Ncppp}
\int_0^T
  \sqrt{\frac{\mathscr{N}(t,\xi)}
             {\mathscr{D}(t,\xi)}\,} \, dt
  \logestimate \, ,
\end{equation}
where $\mathscr{N}$ can be any of~the symbol
\begin{align*}
&  \bigl|\widecheck{N}\bigl(t,\sigma_j(t,\xi),\xi\bigr)\bigr| \, ,
   \quad \text{for $j=1,2$} \, ,  &
&  \bigl|\widecheck{N}\bigl(t,\mu_j(t,\xi),\xi\bigr)\bigr| \, ,
   \quad \text{for $j=1,2$} \, ,  \\
&  \bigl|\widecheck{N}\bigl(t,\tau_j(t,\xi),\xi\bigr)\bigr| \, ,
   \quad \text{for $j=1,2,3$} \, , &
&  \bigl|\widecheck{N}\bigl(t,\lambda_j(t,\xi),\xi\bigr)\bigr| \, ,
   \quad \text{for $j=1,2,3$} \, ,
\end{align*}
and $\mathscr{D}$ can be any of~the symbol
\begin{align*}
&  \bigl|\sigma_2(t,\xi)-\sigma_1(t,\xi)\bigr|+1  &
   &  \bigl|\mu_2(t,\xi)-\mu_1(t,\xi)\bigr|  \\
&  \sqrt{\Delta_L^{(1)}(t,\xi) \,} + 1  &
   &  \sqrt{\Delta_{\mathcal{L}}^{(1)}(t,\xi) \,}  \\
&  \bigl|\tau_3(t,\xi)-\tau_1(t,\xi)\bigr|+1  &
   &  \bigl|\lambda_3(t,\xi)-\lambda_1(t,\xi)\bigr| \, .
\end{align*}

Indeed, from~Lemma~\ref{L-A-1} we have
\[
(\tau_{1}-\tau_{2})^2+(\tau_{2}-\tau_{3})^2+(\tau_{3}-\tau_{1})^2
  =  \frac{9}{2} \, (\sigma_{2}-\sigma_{1})^2 \, ,
\]
hence
\[
|\sigma_2-\sigma_1|+1
  \approx  \sqrt{\Delta_L^{(1)}\,} + 1 \, .
\]

Next, from the elementary inequality
\[
\max(\alpha,\beta,\gamma)
  \le  \sqrt{\alpha^2+\beta^2+\gamma^2\,}
  \le  \sqrt{3\,} \, \max(\alpha,\beta,\gamma) \, ,
\qquad
\text{for any $\alpha,\beta,\gamma\ge0$} \, ,
\]
we see that
\begin{equation} \label{E-91}
|\tau_{3}-\tau_{1}|
  \le  \sqrt{(\tau_{3}-\tau_{1})^2+(\tau_{3}-\tau_{2})^2+(\tau_{2}-\tau_{1})^2\,}
  \le  \sqrt{3\,} \, |\tau_{3}-\tau_{1}| \, ,
\end{equation}
and, consequently
\[
\sqrt{\Delta_L^{(1)}\,} + 1
  \approx  |\tau_{3}-\tau_{1}| + 1 \, .
\]

We prove that
\[
|\mu_2-\mu_1|
   \approx  \sqrt{\Delta_{\mathcal{L}}^{(1)}\,}
   \approx  |\lambda_3-\lambda_1|
\]
in a similar way.

Since
\[
\Bigl| \bigl|\widecheck{N}(\sigma_1)\bigr| - \bigl|\widecheck{N}(\tau_1)\bigr|\Bigr|
  \lesssim  |\sigma_1-\tau_1| \, ,
\]
and
\[
|\sigma_1-\tau_1|
  \le  |\tau_3-\tau_1| \, ,
\]
we see that
\[
\frac{\bigl|\widecheck{N}(\sigma_1)\bigr|}{|\tau_3-\tau_1|+1}
  \approx  \frac{\bigl|\widecheck{N}(\tau_1)\bigr|}{|\tau_3-\tau_1|+1} + 1 \, .
\]

We can prove the other equivalences in a similar way.
\end{Remark}

\begin{Remark} \label{R-M}
By similar arguments we can prove that \eqref{E-dMct} is equivalent to
\begin{equation} \label{E-dMctt}
\int_0^T \frac{\bigl|\partial_\tau \widecheck{M}\bigl(t,\sigma_1(t,\xi),\xi\bigr)\bigr|
               + \bigl|\partial_\tau \widecheck{M}\bigl(t,\sigma_2(t,\xi),\xi\bigr)\bigr|}
              {\bigl|\sigma_1(t,\xi)-\sigma_2(t,\xi)\bigr|+1} \, dt
  \lesssim \log\bigl(1+|\xi|\bigr) \, ,
\end{equation}
where~$\sigma_1$ and $\sigma_2$ are the roots
of~$\partial_\tau L(t,\tau,\xi)$.

More generally,
Hypothesis~\eqref{E-dMct} is equivalent to each of~the conditions
\[
\int_0^T
  \sqrt{\frac{\mathsf{M}(t,\xi)}
             {\mathscr{D}(t,\xi)}\,} \, dt
  \logestimate \, ,
\]
where $\mathsf{M}$ can be any of~the symbol
\begin{align*}
&  \bigl|\partial_\tau \widecheck{M}\bigl(t,\sigma_j(t,\xi),\xi\bigr)\bigr| \, ,
   \quad \text{for $j=1,2$} \, ,  &
&  \bigl|\partial_\tau \widecheck{M}\bigl(t,\mu_j(t,\xi),\xi\bigr)\bigr| \, ,
   \quad \text{for $j=1,2$} \, ,  \\
&  \bigl|\partial_\tau \widecheck{M}\bigl(t,\tau_j(t,\xi),\xi\bigr)\bigr| \, ,
   \quad \text{for $j=1,2,3$} \, , &
&  \bigl|\partial_\tau \widecheck{M}\bigl(t,\lambda_j(t,\xi),\xi\bigr)\bigr| \, ,
   \quad \text{for $j=1,2,3$} \, ,
\end{align*}
and $\mathscr{D}$ is as in the previous Remark.
\end{Remark}

\section{The case of~analytic coefficients} \label{S-an}

In this section we show that if the coefficients are analytic,
then hypothesis~\eqref{E-1}, \eqref{E-2}, \eqref{E-3} and~\eqref{E-4}
are satisfied.

\begin{Lemma}[\cite{CJS}, \cite{Orru1997}, \cite{JannelliT2014}] \label{L-4.1}
Let $f_1,\dotsb,f_d$ be analytic functions
on an open set $\mathscr{O}\subset\C$, $f_j\not\equiv0$ for $j = 1,\dotsb,d$.
For any $\alpha = (\alpha_1,\dotsb,\alpha_d)\in\R^d$ set
\[
\varphi_\alpha(x) := \sum_{j=1}^d \alpha_j f_j(x) \, .
\]

Then for any compact set $\mathscr{K} \subset \mathscr{O}$
there exists $\nu\in\N$ such that either $\varphi_\alpha(x) \equiv 0$
or~$\varphi_\alpha(x)$ has at most $\nu$ zeros in~$\mathscr{K}$,
if counted with their multiplicity.
\end{Lemma}

\begin{proof}
With no loss of~generality we can assume that
$f_1,\dotsb,f_d$ are linearly independent.

If, for any $k\in\N$, there exists $\alpha^{(k)}$ with $\|\alpha^{(k)}\|=1$ such that
$\varphi_\alpha(x)$ has at least $k$ zeros in~$\mathscr{K}$,
by passing to a suitable subsequence,
we can assume that $\alpha^{(k)}$ converges to some $\alpha^*$.
Hence $\varphi_{\alpha^*}$ must have an infinite number of~zeros and hence
it is identically zero.
This contradicts the fact that the $f_1,\dotsb,f_d$ are linearly independent.
\end{proof}

As any polynomial $\Phi(t,\xi)$ in~$\xi$ with analytic coefficients
can be regarded as a linear combination of~its coefficients
we deduce from Lemma~\ref{L-4.1}

\begin{Corollary} \label{C-4.2}
Let $\Phi(t,\xi)$ be a polynomial in~$\xi$ with analytic coefficients
on an open set $\mathscr{O}\subset\C$.

Then for any compact set $\mathscr{K} \subset \mathscr{O}$
there exists $\nu\in\N$ such that either $\Phi(t,\xi) \equiv 0$
or~$\Phi(\cdot,\xi)$ has at most $\nu$ zeros in~$\mathscr{K}$,
if counted with their multiplicity.
\end{Corollary}

\begin{Proposition} \label{P-4.3}
Let $\PP(t,\tau,\xi)$ be a third order monic hyperbolic polynomial in~$\tau$,
whose coefficients are polynomial in~$\xi$ and analytic in~$t\in\mathscr{O}$,
$\mathscr{O}$ open set in~$\C$.
Assume that the roots in~$\tau$ of~$\PP(t,\tau,\xi)$ are distinct for $\xi$ large:
\[
\tau_1(t,\xi) < \tau_2(t,\xi) < \tau_3(t,\xi) \, ,
\qquad
\text{if $|\xi|\ge R$} \, .
\]

Let $\QQ(t,\tau,\xi) $ be another polynomial in~$\tau$
whose coefficients are polynomial in~$\xi$ and analytic in~$t\in\mathscr{O}$.

Then for any compact set $\mathscr{K} \subset \mathscr{O}$
there exists $\nu\in\N$ such that for any $\xi\in\R^n$,
with $|\xi|\ge R$,
the functions
\[
t \mapsto \QQ\bigl(t,\tau_j(t,\xi),\xi\bigr) \, ,
\qquad
j=1,2,3 \, ,
\]
are either identically zero or~have at most $\nu$ zeros in~$\mathscr{K}$.
\end{Proposition}

\begin{Remark}
If a root $\tau_j(t,\xi)$ were a linear function of~$\xi$,
then $\QQ\bigl(t,\tau_j(t,\xi),\xi\bigr)$ would be a polynomial in~$\xi$,
and, by Corollary~\ref{C-4.2}, we can get easly the~result.

Unfortunately, in general, the roots $\tau_j(t,\xi)$ are not linear function of~$\xi$,
hence we cannot apply Corollary~\ref{C-4.2} directly.
\end{Remark}

\begin{proof}
First of~all we remark that, by the implicit function theorem,
the roots $\tau_j=\tau_j(t,\xi)$ are analytic functions for $\xi$ large,
thus so are the functions $\QQ\bigl(t,\tau_j(t,\xi),\xi\bigr)$, $j=1,2,3$.

Let
\begin{align*}
\Phi_1(t,\xi)
&  =  \QQ\bigl(t,\tau_1(t,\xi),\xi\bigr)+\QQ\bigl(t,\tau_2(t,\xi),\xi\bigr)+\QQ\bigl(t,\tau_3(t,\xi),\xi\bigr)  \\
\Phi_2(t,\xi)
&  =  \QQ\bigl(t,\tau_1(t,\xi),\xi\bigr) \, \QQ\bigl(t,\tau_2(t,\xi),\xi\bigr)  \\
&  \qquad  +  \QQ\bigl(t,\tau_2(t,\xi),\xi\bigr) \, \QQ\bigl(t,\tau_3(t,\xi),\xi\bigr)
      +  \QQ\bigl(t,\tau_3(t,\xi),\xi\bigr) \, \QQ\bigl(t,\tau_1(t,\xi),\xi\bigr)  \\
\Phi_3(t,\xi)
&  =  \QQ\bigl(t,\tau_1(t,\xi),\xi\bigr) \, \QQ\bigl(t,\tau_2(t,\xi),\xi\bigr) \,
      \QQ\bigl(t,\tau_3(t,\xi),\xi\bigr) \, .
\end{align*}

These functions are symmetric in the~$\tau_j$'s
and so they are polynomials in~$\xi$,
with analytic coefficients in~$t$.

For a given compact $\mathscr{K} \subset \mathscr{O}$,
consider another compact $\mathscr{L} \subset \mathscr{O}$
such that $\mathscr{K} \subset \ring{\mathscr{L}}$.
By~Co\-rollary~\ref{C-4.2},
there exists $\nu\in\N$ such that, for any~$j=1,2,3$ and $\xi\in\R^n$
either $\Phi_j(\cdot,\xi)$ is identically zero
or $\Phi_j(\cdot,\xi)$ has at most~$\nu$ zeros in~$\mathscr{L}$.

Assume at first that $\Phi_3(\cdot,\xi) \not\equiv 0$.

Let $\xi_0\in\R^n\setminus\{0\}$ be fixed,
if $\Phi_3(\cdot,\xi_0) \not\equiv 0$,
then $\Phi_3(\cdot,\xi_0)$ has at most $\nu$ zeros in~$\mathscr{L}$,
hence, also the $\QQ\bigl(t,\tau_j(t,\xi_0),\xi_0\bigr)$'s have at most $\nu$ zeros in~$\mathscr{L}$.
If $\Phi_3(\cdot,\xi_0) \equiv 0$
but one of~the $\QQ\bigl(t,\tau_j(t,\xi_0),\xi_0\bigr)$ is not identically zero,
then it has $\eta$ zeros in~$\mathscr{L}$.
Since $\Phi_3(\cdot,\xi) \not\equiv 0$,
we can find a sequence $\{\xi^{(k)}\}_{k\in\N}$
convergent to $\xi_0$
and such that $\Phi_3(\cdot,\xi^{(k)}) \not\equiv 0$.
If $k$ is large enough,
by Rouch\'{e} Theorem,
$\QQ\bigl(t,\tau_j(t,\xi^{(k)}),\xi^{(k)}\bigr)$
has at most~$\nu$ zeros in~$\mathscr{L}$,
then $\QQ\bigl(t,\tau_j(t,\xi_0),\xi_0\bigr)$ has $\eta\le\nu$ zeros in~$\mathscr{K}$.

If $\Phi_3(t,\xi) \equiv 0$,
then at least one of~the $\QQ\bigl(t,\tau_j(t,\xi),\xi\bigr)$
(say $\QQ\bigl(t,\tau_3(t,\xi),\xi\bigr)$) is identically zero,
thus $\Phi_2(\cdot,\xi)$ reduces to
$\QQ\bigl(t,\tau_1(t,\xi),\xi\bigr) \, \QQ\bigl(t,\tau_2(t,\xi),\xi\bigr)$.

We can repeat the same argument as before:
if $\Phi_2(\cdot,\xi)$ is not identical zero,
then $\QQ\bigl(t,\tau_1(t,\xi),\xi\bigr)$ and $\QQ\bigl(t,\tau_2(t,\xi),\xi\bigr)$
have at most $\nu$ zeros in~$\mathscr{K}$,
whereas if $\Phi_2(\cdot,\xi)$ is identical zero,
at least one between $\QQ\bigl(t,\tau_1(t,\xi),\xi\bigr)$ and $\QQ\bigl(t,\tau_2(t,\xi),\xi\bigr)$
vanishes identically.
If $\QQ\bigl(t,\tau_2(t,\xi),\xi\bigr)$ is identically zero,
then $\Phi_1(t,\xi) = \QQ\bigl(t,\tau_1(t,\xi),\xi\bigr)$,
and we get immediately that $\QQ\bigl(t,\tau_1(t,\xi),\xi\bigr)$ is
either identically zero or it has at most $\nu$ zeros in~$\mathscr{K}$.
\end{proof}

\begin{Proposition} \label{P-4.4}
Let $\PP(t,\tau,\xi)$ be as in Proposition~\ref{P-4.3}
and let $\QQ(t,\tau,\sigma,\xi)$ be a symmetric polynomial in~$(\tau,\sigma)$
whose coefficients are polynomial in~$\xi$ and analytic in~$t\in\mathscr{O}$.

Then for any compact set $\mathscr{K} \subset \mathscr{O}$
there exists $\nu\in\N$ such that for any $\xi$
the functions
\[
t \mapsto \QQ\bigl(t,\tau_j(t,\xi),\tau_k(t,\xi),\xi\bigr) \, ,
\qquad
(j,k)\in\SS_2 \, ,
\]
are either identically zero or~have at most $\nu$ zeros in~$\mathscr{K}$.
\end{Proposition}

\goodbreak

The proof is similar to that of~Proposition~\ref{P-4.3},
and is obtained by considering the functions
\begin{align*}
\Phi_1(t,\xi)
&  =  \QQ_{1,2}(t,\xi)+\QQ_{2,3}(t,\xi)+\QQ_{3,1}(t,\xi)  \\
\Phi_2(t,\xi)
&  =  \QQ_{1,2}(t,\xi) \, \QQ_{2,3}(t,\xi) + \QQ_{2,3}(t,\xi) \, \QQ_{3,1}(t,\xi)
      +  \QQ_{3,1}(t,\xi) \, \QQ_{1,2}(t,\xi)  \\
\Phi_3(t,\xi)
&  =  \QQ_{1,2}(t,\xi) \, \QQ_{2,3}(t,\xi) \, \QQ_{3,1}(t,\xi) \, ,
\end{align*}
where, for the sake of~brevity, we have set
\[
\QQ_{j,k}(t,\xi)
  =  \QQ\bigl(t,\tau_j(t,\xi),\tau_k(t,\xi),\xi\bigr) \, .
\]

\begin{Proposition} \label{P-1}
Let
\[
\Xi
  =  \Bigl\{ \ \xi \in \R^n \ \Bigm| \ |\xi| \ge R \ \Bigr\} \, ,
\]
for some $R>0$ and
let $f\colon[0,T]\times\Xi\to\R$ be such that
\begin{enumerate}
\item  $f$ is Lipschitz continuous in $t$, uniformly with respect to $\xi$,
      that is there exists $C_0>0$ such that
         \[
         \bigl|f(t_1,\xi)-f(t_2,\xi)\bigr|
            \le C_0 |t_1-t_2| \, ,
         \qquad
         \text{for any $t_1,t_2\in[0,T]$ and $\xi\in\Xi$} \, ;
         \]

\item  there exist positive constants $A,C_1,C_2$ such that
         \[
         C_1 \, \bigl(1+|\xi|\bigr)^{-A}
           \le  f(t,\xi)
           \le  C_2 \, \bigl(1+|\xi|\bigr)^A
         \]
         for any $t\in[0,T]$ and $\xi\in\Xi$;

\item  there exists $\nu\in\N$ such that
         for any $\xi\in\Xi$ there exists a partition
         $0=t_1<t_2<\dotsb<t_{\mu-1}<t_\mu=T$ of~$[0,T]$,
         with~$\mu\le\nu$ such that
         \begin{itemize}
            \item  $f(\cdot,\xi) \in \mathcal{C}^1\bigl(\left]t_j,t_{j+1}\right[\bigr)$,
                   for $j=1,2,\dotsc,\mu-1$;
            \item  $\partial_t f(t,\xi)\ne0$ for any $t \in \left]t_j,t_{j+1}\right[$,
                   for $j=1,2,\dotsc,\mu-1$.
         \end{itemize}
\end{enumerate}

Then $f$ satisfies the logarithmic condition
\[
\int_0^T \frac{\bigl|\partial_t f(t,\xi)\bigr|}{f(t,\xi)} \, dt
  \logestimate \, ,
\quad
\text{for any $\xi\in\Xi$} \, .
\]
\end{Proposition}

\begin{proof}
Let us fix $\xi$.
As $\partial_t f(t,\xi)$ does not change sign in~$\left]t_j,t_{j+1}\right[$,
we have
\[
\int_{t_j}^{t_{j+1}}
  \frac{\bigl|\partial_t f(t,\xi)\bigr|}{f(t,\xi)} \, dt
  =  \biggl|\int_{t_j}^{t_{j+1}}
            \frac{\partial_t f(t,\xi)}{f(t,\xi)} \, dt\biggr|
  =  \biggl|\log f(t_{j+1},\xi) - \log f(t_j,\xi)\biggr|
 \le  C^* \, \log\bigl(1+|\xi|\bigr) \, ,
\]
where $C^*$ depends on $A,C_1,C_2$.
Hence
\[
\int_0^T
  \frac{\bigl|\partial_t f(t,\xi)\bigr|}{f(t,\xi)} \, dt
 \le  \nu \, C^* \, \log\bigl(1+|\xi|\bigr) \, ,
\]
where $\nu$ and $C^*$ are independent of~$\xi$.
\end{proof}

\subsection*{Study of the condition~\eqref{E-1}}
\text{}

Let
\[
f(t,\xi) = \bigl|\lambda_j (t,\xi) - \lambda_h (t,\xi)\bigr| \, .
\]

The function $f(t,\xi)$ never vanishes and
its critical points verify
\[
\partial_t \bigl(\lambda_j (t,\xi) - \lambda_h (t,\xi)\bigr) = 0 \, .
\]

Now, by the implicit function Theorem,
we~have
\begin{equation} \label{E-dt}
\partial_t \lambda_{j}(t,\xi)
  =  - \frac{(\partial_t L_\e) (\lambda_{j})}
            {(\partial_\tau L_\e) (\lambda_{j})} \, ,
\end{equation}
where, for the sake of~brevity,
we~write
\begin{align*}
(\partial_t L_\e) (\lambda_{j})
&  =  \partial_t L_\e(t,\tau,\xi)\bigm|_{\tau = \lambda_{j}(t,\xi)} \, ,  \\
(\partial_\tau L_\e) (\lambda_{j})
&  =  \partial_\tau L_\e(t,\tau,\xi)\bigm|_{\tau = \lambda_{j}(t,\xi)} \, .
\end{align*}

{}From~\eqref{E-dt} we get
\[
\Bigl[\partial_t\bigl(\lambda_{j}(t_1,\xi) - \lambda_{h}(t_1,\xi)\bigr)\Bigr]^2
  =  \frac{Q(x,\lambda_j,\lambda_h,\xi)}
          {\bigl[(\partial_\tau L_\e) (\lambda_{j})\bigr]^2 \,
           \bigl[(\partial_\tau L_\e) (\lambda_{h})\bigr]^2} \, ,
\]
where
\[
Q(x,\tau,\sigma,\xi)
  =  \Bigl[(\partial_t L_\e) (\tau) \,
           (\partial_\tau L_\e) (\sigma)  \\
           -  (\partial_t L_\e) (\sigma) \,
              (\partial_\tau L_\e) (\tau)\Bigr]^2 \, .
\]

The polynomial $Q$ verifies the hypothesis of~Proposition~\ref{P-4.4},
hence the number of~zeros of~the function $t\mapsto\partial_t f(t,\xi)$
is bounded w.r.t.~$\xi\in\Xi$,
and, applying Proposition~\ref{P-1} to $f$,
we see that~\eqref{E-1} holds true.

\subsection*{Study of the condition~\eqref{E-2}}
\text{}

Let
\[
f(t,\xi)
  =  \bigl|\partial_t\lambda_{j}(t_1,\xi) - \partial_t\lambda_{h}(t_1,\xi)\bigr| + 1 \, .
\]
If $\partial_t f(t,\xi)$ changes sign at~$t_1$, then
either
$\partial_t\lambda_{j}(t_1,\xi) - \partial_t\lambda_{h}(t_1,\xi)=0$
or
$\partial_t^2\lambda_{j}(t_1,\xi) - \partial_t^2\lambda_{h}(t_1,\xi)=0$.

The first case can be treated as before.
For the second, from~\eqref{E-dt} we get
\[
\partial_{tt}^2 \lambda_{j}(t,\xi)
  =  \frac{\psi(t,\lambda_j,\xi)}{\bigl[(\partial_\tau L_\e) (\lambda_{j})\bigr]^3} \, ,
\]
where
\begin{align*}
\psi(t,\tau,\xi)
&  =  2 \partial_{t\tau}^2 L_\e (\tau) \,
       \partial_t L_\e (\tau) \,
       \partial_\tau L_\e (\tau)  \\
&  \qquad
      -  \partial_{\tau\tau}^2 L_\e (\tau) \, \bigl[\partial_t L_\e (\tau)\bigr]^2
      -  \partial_{tt}^2 L_\e (\tau) \, \bigl[\partial_\tau L_\e (\tau)\bigr]^2 \, ,
\end{align*}
hence
\[
\Bigl[\partial_{tt}^2\bigl(\lambda_{j}(t_1,\xi) - \lambda_{h}(t_1,\xi)\bigr)\Bigr]^2
  =  \frac{Q(x,\lambda_{j},\lambda_{h},\xi)}
          {\bigl[(\partial_\tau L_\e) (\lambda_{j})\bigr]^6 \,
           \bigl[(\partial_\tau L_\e) (\lambda_{h})\bigr]^6} \, ,
\]
where
\[
Q(x,\tau,\sigma,\xi)
  =  \Bigl[\psi(t,\tau,\xi) \, \bigl[\partial_\tau L_\e (\sigma)\bigr]^3  \\
           -  \psi(t,\sigma,\xi) \, \bigl[\partial_\tau L_\e (\tau)\bigr]^3\Bigr]^2 \, .
\]

The polynomial $Q$ verifies the hypothesis of~Proposition~\ref{P-4.4},
hence the number of~zeros of~the function $t\mapsto\partial_t f(t,\xi)$
is bounded w.r.t.~$\xi\in\Xi$,
and, applying Proposition~\ref{P-1} to $f$,
we see that~\eqref{E-2} holds true.

\subsection*{Study of the condition~\eqref{E-3}}
\text{}

If $M$ has real coefficients,
we consider
\[
f(t,\xi)
  =  \bigl|\widecheck{M}\bigl(t,\lambda_{j}(t,\xi),\xi\bigr)\bigl|+1 \, .
\]
For fixed~$\xi$,
the oscillations of~$f(t,\xi)$ are zeros in~$t$ of~either
$\widecheck{M}\bigl(t,\lambda_{j}(t,\xi),\xi\bigr)$
or
$\partial_t \widecheck{M}\bigl(t,\lambda_{j}(t,\xi),\xi\bigr)$.

Now
\[
\partial_t \widecheck{M}\bigl(t,\lambda_{j}(t,\xi),\xi\bigr)
  =  \partial_t \widecheck{M}(t,\tau,\xi)\bigm|_{\tau = \lambda_{j}(t,\xi)}
     +  \partial_\tau \widecheck{M}(t,\tau,\xi)\bigm|_{\tau = \lambda_{j}(t,\xi)}
        \partial_t \lambda_{j}(t,\xi) \, ,
\]
and,
by~\eqref{E-dt} we see that
$\partial_t \widecheck{M}\bigl(t,\lambda_{j}(t,\xi),\xi\bigr) = 0$
if, and only if, $Q\bigl(t,\lambda_{j}(t,\xi),\xi\bigr)=0$, where
\[
Q(t,\tau,\xi)
  =  \partial_t \widecheck{M}(t,\tau,\xi)
     \partial_\tau L_\e(t,\tau,\xi)
     -  \partial_\tau \widecheck{M}(t,\tau,\xi)
        \partial_t L_\e(t,\tau,\xi) \, .
\]

The polynomial $Q$ verifies the hypothesis of~Proposition~\ref{P-4.4},
hence the number of~zeros of~the function $t\mapsto\partial_t f(t,\xi)$
is bounded w.r.t.~$\xi\in\Xi$,
and, applying Proposition~\ref{P-1} to $f$,
we see that~\eqref{E-3} holds true.

\smallskip

If the coefficients of~$\widecheck{M}$ are complex
we~consider the zeros of
\[
\partial_t |\widecheck{M}|^2
  =  2\Re \bigl(\partial_t \widecheck{M} \, \overline{\widecheck{M}}\bigr)
  =  2 \partial_t \Re(\widecheck{M}) \, \Re(\widecheck{M})
     +  2 \partial_t \Im(\widecheck{M}) \, \Im(\widecheck{M})
\]
and by the same argument we get~\eqref{E-3} again.

The proof that condition~\eqref{E-4} is satisfied,
is similar to that of~condition~\eqref{E-3},
so we omit~it.

\section{Pointwise Levi conditions} \label{S-anLevi}

Throughout this section we assume that the coefficients
of the operator are analytic,
and we~express the Levi conditions~\eqref{E-Mcl} and~\eqref{E-Ncl}
as pointwise conditions.

\smallskip

We have to distinguish three cases:
\begin{description}
\item[Case~I]    $\Delta_L \not \equiv 0$;
\item[Case~II]   $\Delta_L \equiv 0$ and $\Delta_L^{(1)} \not\equiv0$;
\item[Case~III]  $\Delta_L \equiv \Delta_L^{(1)}\equiv 0$.
\end{description}

\subsection*{Case~I: $\Delta_L \protect\not \equiv 0$}
\text{}

We consider at first the terms of~order 2.

\begin{Proposition} \label{P-Ma}
Assume that
\begin{equation} \label{E-170D}
\bigl|\tau_k(t,\xi) - \tau_l(t,\xi)\bigr| \, \Bigl|\widecheck{M}\bigl(t,\tau_j(t,\xi),\xi\bigr)\Bigr|
  \lesssim  \sqrt{\Delta_L(t,\xi) \,}
            +  \Bigl|\partial_t \sqrt{\Delta_L(t,\xi) \,}\Bigr| \, ,
\end{equation}
for all $(j,k,l)\in\SS_3$,
then condition~\eqref{E-Mcl} is verified.
\end{Proposition}

For the proof we need the following Lemma
\cite[Proposition~4.1]{JannelliT2011}.

\begin{Lemma} \label{L-Prato 2}
Let $\Delta(t,\xi)$ be an homogeneous polynomial in~$\xi$ with coefficient analytic in~$t$
and assume that $\Delta(t,\xi)\not\equiv0$.
Then:
\begin{enumerate}
\item  there exists
       $X\subset \mathbb{S}^n := \bigl\{ \, \xi\in\R^n \, \bigm| \, |\xi|=1 \, \bigr\}$
       such that $\Delta(t,\xi)\not\equiv0$ in $\left]-\delta,T+\delta\right[$ for any $\xi\in X$,
       and the set $\mathbb{S}^n\setminus X$ is negligeable with respect
       to the Hausdorff $(n-1)$--measure;

\item  for any~$[a,b]\subset\left]-\delta,T+\delta\right[$
       we~can find constants~$c_1,c_2>0$ and~$p,q\in\N$
       such that for~any~$\xi\in X$ and any~$\e\in(0,1/e]$
       there exists~$A_{\xi,\e}\subset[a,b]$ such that
       \begin{enumerate}
       \item  $A_{\xi,\e}$ is a union of~at most~$p$ disjoint intervals,
       \item  $m(A_{\xi,\e}) \le \e$,
       \item  $\min_{t \in [a,b]\setminus A_{\xi,\e}} \Delta(t,\xi)
                 \ge  c_1 \, \e^{2q}
                      \bigl\|\Delta(\cdot,\xi)\bigr\|_{L^\infty([a,b])}$
       \item  $\int\limits_{[a,b]\setminus A_{\xi,\e}}
                        \frac{\bigl|\Delta'(t,\xi)\bigr|}{\Delta(t,\xi)}\,dt
                   \le  c_2 \, \log\frac{1}{\e}$.
       \end{enumerate}
\end{enumerate}
\end{Lemma}

\begin{proof}[Proof of~Proposition~\ref{P-Ma}]
We will prove that~\eqref{E-170D} implies~\eqref{E-Mct} and~\eqref{E-dMct}
hence, thanks to Proposition~\ref{P-M}, we get~\eqref{E-Mcl}.

Let $\e=|\xi|^{-2}$,
with $|\xi|$ large enough,
and let $A_{\xi,\e}$ be the set given by Lemma~\ref{L-Prato 2}
with $\Delta(t,\xi) = \Delta_L(t,\xi)$;
we have
\begin{equation} \label{E-le e}
\int_{A_{\xi,\e}} \frac{\bigl|\widecheck{M}(\tau_j)\bigr|}
                  {\bigl(|\tau_j-\tau_k|+1\bigr)
                   \bigl(|\tau_j-\tau_l|+1\bigr)} \, dt
  \lesssim  \int_{A_{\xi,\e}} |\xi|^2 \, dt
  \lesssim  \e\,|\xi|^2
   = 1 \, .
\end{equation}

On the other side,
as
\[
\frac{\bigl|\widecheck{M}(\tau_j)\bigr|}
     {\bigl(|\tau_j-\tau_k|+1\bigr)
      \bigl(|\tau_j-\tau_l|+1\bigr)}
  \le  \frac{|\tau_k-\tau_l| \,
             \bigl|\widecheck{M}(\tau_j)\bigr|}
            {\sqrt{\Delta_L(t,\xi) \,}} \, ,
\]
assuming~\eqref{E-170D},
thanks to Lemma~\ref{L-Prato 2},
we have
\begin{equation} \label{E-ge e}
\begin{split}
\int_{\left]-\delta,T+\delta\right[ \setminus A_{\xi,\e}}
&    \frac{\bigl|\widecheck{M}(\tau_j)\bigr|}
          {\bigl(|\tau_j-\tau_k|+1\bigr)
           \bigl(|\tau_j-\tau_l|+1\bigr)} \, dt  \\
&  \lesssim  \int_{\left]-\delta,T+\delta\right[ \setminus A_{\xi,\e}}
                  1 + \frac{\bigl|\partial_t \Delta_L(t,\xi)\bigr|}
                           {\Delta_L(t,\xi)} \, dt
   \logestimate \, .
\end{split}
\end{equation}

Combining~\eqref{E-le e} and~\eqref{E-ge e} we get~\eqref{E-Mct}.

Now we prove that~\eqref{E-dMct} holds true.
As the roots of~$L$ are distinct for a.e.~$(t,\xi)$,
by the Lagrange interpolation formula
we~have
\begin{equation} \label{E-845}
\widecheck{M}(t,\tau,\xi)
  =  \ell_1(t,\xi) \, L_{23}(t,\tau,\xi)
     +  \ell_2(t,\xi) \, L_{13}(t,\tau,\xi)
     +  \ell_3(t,\xi) \, L_{12}(t,\tau,\xi) \, ,
\end{equation}
where the operators $L_{jk}$ are the operators $L_{jk,\e}$
defined in~\textsection\ref{S-proof} with $\e=0$,
the $\ell_j$ are given by
\begin{equation} \label{E-ell}
\ell_j(t,\xi)
  :=  \frac{\widecheck{M}\bigl(t,\tau_j(t,\xi),\xi\bigr)}
     {\bigl(\tau_j(t,\xi)-\tau_h(t,\xi)\bigr)
      \bigl(\tau_j(t,\xi)-\tau_l(t,\xi)\bigr)}
\end{equation}
hence, by Hypothesis~\eqref{E-170D},
they verify
\begin{equation} \label{E-601}
\bigl|\ell_j(t,\xi)\bigr|
  \lesssim  1 + \frac{\bigl|\partial_t \Delta_L(t,\xi)\bigr|}
                     {\Delta_L(t,\xi)} \, .
\end{equation}

On the other side, differentiating~\eqref{E-845} w.r.t~$\tau$,
we~get
\begin{align*}
\partial_\tau \widecheck{M}(\tau)
&  =  \ell_1 \, (L_2+L_3)
      +  \ell_2 \, (L_3+L_1)
      +  \ell_3 \, (L_1+L_2)  \\
&  =  (\ell_2+\ell_3) \, L_1
      +  (\ell_3+\ell_1) \, L_2
      +  (\ell_1+\ell_2) \, L_3 \, ,
\end{align*}
where the operators $L_j$ are the operators $L_{j,\e}$
defined in~\textsection\ref{S-proof} with $\e=0$,

Now, since $\tau_1 \le \tau_2 \le \tau_3$,
we can find $\theta\in[0,1]$ such that
\[
\tau_2
  =  \theta \, \tau_1 + (1 - \theta) \, \tau_3 \, ,
\]
hence
\[
L_2
  =  \theta \, L_1
     +  (1 - \theta) \, L_3 \, ,
\]
consequently we can write
\begin{equation} \label{E-751}
\partial_\tau \widecheck{M}(\tau)
  =  \wt{\ell}_3 \, L_1(\tau)
      +  \wt{\ell}_1 \, L_3(\tau) \, ,
\end{equation}
where $\wt{\ell}_1$ and $\wt{\ell}_3$
are some linear combination with bounded coefficients
of $\ell_1$, $\ell_2$ and~$\ell_3$,
hence the $\wt{\ell}_j$s verify~\eqref{E-601} too.

Let $\e=|\xi|^{-1}$,
with $|\xi|$ large enough,
and let $A_{\xi,\e}$ be as above,
we have
\begin{align*}
\int_{A_{\xi,\e}} \frac{\bigl|\partial_\tau \widecheck{M}(\tau_1)\bigr|
                   +  \bigl|\partial_\tau \widecheck{M}(\tau_3)\bigr|}
                 {|\tau_1-\tau_3|+1} \, dt
&  \lesssim  \int_{A_{\xi,\e}} |\xi| \, dt
   \lesssim  \e\,|\xi|
   = 1 \, ,
\intertext{%
whereas, thanks to~\eqref{E-751}:}
\int_{\left]-\delta,T+\delta\right[ \setminus A_{\xi,\e}}
             \frac{\bigl|\partial_\tau \widecheck{M}(\tau_1)\bigr|
                    +  \bigl|\partial_\tau \widecheck{M}(\tau_3)\bigr|}
                  {|\tau_1-\tau_3|+1} \, dt
&  \lesssim  \int_{\left]-\delta,T+\delta\right[ \setminus A_{\xi,\e}}
                  \bigl|\wt{\ell}_1(t,\xi)\bigr|
                  +  \bigl|\wt{\ell}_3(t,\xi)\bigr| \, dt  \\
&  \lesssim  \int_{\left]-\delta,T+\delta\right[ \setminus A_{\xi,\e}}
                  1 + \frac{\bigl|\partial_t \Delta_L(t,\xi)\bigr|}
                           {\Delta_L(t,\xi)} \, dt  \\
&  \logestimate \, .
\end{align*}

Combining the above estimates we get~\eqref{E-dMct}.
\end{proof}

Condition~\eqref{E-170D} means that if, for fixed~$\xi$,
the function
\[
t \mapsto \bigl(\tau_j(t,\xi) - \tau_k(t,\xi)\bigr) \,
          \bigl(\tau_j(t,\xi) - \tau_l(t,\xi)\bigr)
\]
vanishes of order~$\nu$ at~$t=t_0$,
then the function
\[
t \mapsto \widecheck{M}\bigl(t,\tau_j(t,\xi),\xi\bigr)
\]
must vanish (at least) at order $\nu-1$ at~$t=t_0$.
Thus in space dimension~$n=1$ Proposition~\ref{P-Ma} can be precised.

\begin{Proposition} \label{P-Ma1}
In space dimension $n=1$,
\eqref{E-Mcl} is equivalent to the following condition:

there exist $t_1,\dotsc,t_\nu \in \left]-\delta,T+\delta\right[$
such that
\begin{equation} \label{E-170}
\prod_{h=1}^\nu |t-t_h| \,
\Bigl|\widecheck{M}\bigl(t,\tau_j(t,\xi),\xi\bigr)\Bigr|
  \lesssim  \bigl|\tau_j(t,\xi)-\tau_k(t,\xi)\bigr| \,
            \bigl|\tau_j(t,\xi)-\tau_l(t,\xi)\bigr| \, ,
\end{equation}
with $j,k,l$ such that $(j,k,l)\in\SS_3$.
\end{Proposition}

\begin{Remark}
In the special case $\nu=1$ and $t_1=0$,
condition~\eqref{E-170}
reduces to
\begin{equation} \label{E-70}
|t|\,\Bigl|\widecheck{M}\bigl(t,\tau_j(t,\xi),\xi\bigr)\Bigr|
  \lesssim  \bigl|\tau_j(t,\xi)-\tau_k(t,\xi)\bigr| \,
            \bigl|\tau_j(t,\xi)-\tau_l(t,\xi)\bigr| \, ,
\end{equation}
for any $(j,k,l)\in\SS_3$.
\end{Remark}

\begin{proof}
By the previous remark,
we see that~\eqref{E-170} is equivalent to~\eqref{E-170D}
and implies~\eqref{E-Mct} and~\eqref{E-dMct}.

Now we prove that~\eqref{E-Mct} implies~\eqref{E-170} by contradiction.
First of~all, as the zeros of~$\Delta_L$ are isolated,
we can decompose $\left]-\delta,T+\delta\right[$
into a finite number of~contiguous subintervals
each containing a zero of~$\Delta_L$.
Thus, with no loss of~generality,
we can restrict to the case in which
$\Delta_L(t)$ vanishes only at~$t=0$;
in this case condition~\eqref{E-170} reduces to~\eqref{E-70}.

Suppose that~\eqref{E-70} is violated,
hence, with no loss of~generality, we have
\[
\frac{\Bigl|\widecheck{M}\bigl(t,\tau_3(t) \,\xi,\xi\bigr)\Bigr|}
     {\bigl|\tau_3(t)-\tau_1(t)\bigl|
      \bigl|\tau_3(t)-\tau_2(t)\bigl|\,|\xi|^2}
  \gtrsim  \frac{1}{t^m}
\]
for some $m\ge 2$.

As $\Delta_L(0)=0$ and $\Delta_L(t)\ne0$ for $t\ne0$,
there exist~$r_1,r_2$ such that
\[
\bigl|\tau_3(t)-\tau_1(t)\bigl| \gtrsim t^{r_1}
\qquad
\bigl|\tau_3(t)-\tau_2(t)\bigl| \gtrsim t^{r_2}
\]
and $r \defeq \max(r_1,r_2) \ge 1$.
Note that $\min(r_1,r_2)>0$
if and only if $t=0$ is a triple point.

For $t\ge\e^{1/r}$, $\e=\frac{1}{|\xi|}$ and $|\xi|\ge1$
we have
\[
\bigl|\tau_3(t)-\tau_1(t)\bigl|\,|\xi| \gtrsim 1
\qquad
\bigl|\tau_3(t)-\tau_2(t)\bigl|\,|\xi| \gtrsim 1 \, ,
\]
hence
\begin{align*}
\int_{\e^{1/r}}^T
&   \frac{\Bigl|\widecheck{M}\bigl(t,\tau_3(t) \,\xi,\xi\bigr)\Bigr|}
         {\Bigl(\bigl|\tau_3(t)-\tau_1(t)\bigl|\,|\xi|+1\Bigr)
          \Bigl(\bigl|\tau_3(t)-\tau_2(t)\bigl|\,|\xi|+1\Bigr)} \, dt  \\
&  \gtrsim
        \int_{\e^{1/r}}^T \frac{1}{t^m} \,
        \frac{\bigl|\tau_3(t)-\tau_1(t)\bigl|\,|\xi|}
             {\bigl|\tau_3(t)-\tau_1(t)\bigl|\,|\xi|+1} \,
        \frac{\bigl|\tau_3(t)-\tau_2(t)\bigl|\,|\xi|}
             {\bigl|\tau_3(t)-\tau_2(t)\bigl|\,|\xi|+1} \, dt
  \gtrsim \int_{\e^{1/r}}^T \frac{1}{t^m} \, dt
  \approx |\xi|^{\frac{m-1}{r}} \, .
\end{align*}
This shows that~\eqref{E-Mct} cannot hold true if~\eqref{E-70} is violated.
\end{proof}

Now we consider the terms of~order 1.

\begin{Proposition} \label{P-Na}
Assume that
\begin{equation} \label{E-NDD}
\Bigl|\widecheck{N}\bigl(t,\sigma_j(t,\xi),\xi\bigr)\Bigr|
  \lesssim  \sqrt{ \Delta_{\partial_t L}(t,\xi) \, }
             + \frac{\bigl(\partial_t \Delta_{\partial_t L}(t,\xi)\bigr)^2}
                    {\bigl[ \Delta_{\partial_t L}(t,\xi) \bigr]^{3/2}} \, ,
\qquad
j=1,2 \, ,
\end{equation}
then condition~\eqref{E-Ncl} is verified.
\end{Proposition}

\begin{proof}
First of~all we recall that, by~Proposition~\ref{P-Ncpp},
\eqref{E-Ncl} is equivalent to~\eqref{E-Ncpp},
hence we will prove that~\eqref{E-NDD} implies~\eqref{E-Ncpp}.

Let $\e=|\xi|^{-1/2}$,
with $|\xi|$ large enough,
and let $A_{\xi,\e}$ be the set given by Lemma~\ref{L-Prato 2}
with $\Delta(t,\xi) = \Delta_{\partial_t L}(t,\xi)$;
we have
\[
\int_{A_{\xi,\e}} \sqrt{\frac{\bigl|\widecheck{N}\bigl(t,\sigma_j(t,\xi),\xi\bigr)\bigr|}
                       {\bigl|\sigma_2(t,\xi) - \sigma_1(t,\xi)\bigr|+1}\,} \, dt
  \lesssim  \int_{A_{\xi,\e}} |\xi|^{1/2} \, dt
  \lesssim  \e\,|\xi|^{1/2}
   = 1 \, .
\]

On the other side,
assuming~\eqref{E-NDD}
and thanks to Lemma~\ref{L-Prato 2},
we have
\begin{align*}
\int_{\left]-\delta,T+\delta\right[ \setminus A_{\xi,\e}}
&             \sqrt{\frac{\bigl|\widecheck{N}\bigl(t,\sigma_j(t,\xi),\xi\bigr)\bigr|}
           {\bigl|\sigma_2(t,\xi) - \sigma_1(t,\xi)\bigr|+1}\,} \, dt  \\
&  \lesssim  \int_{\left]-\delta,T+\delta\right[ \setminus A_{\xi,\e}}
                  1 + \frac{\bigl|\partial_t \Delta_{\partial_t L}(t,\xi)\bigr|}
                           {\Delta_{\partial_t L}(t,\xi)} \, dt
  \logestimate \, .
\end{align*}

Combining the above estimates we get~\eqref{E-Ncpp}.
\end{proof}

Condition~\eqref{E-NDD} means that if, for fixed~$\xi$,
the function $t \mapsto \Delta_{\partial_t L}(t,\xi)$
vanishes of order~$2\nu$, with $\nu>2$, at~$t=t_0$,
then the function $t \mapsto \widecheck{N}\bigl(t,\tau_j(t,\xi),\xi\bigr)$
must vanish (at least) at order $\nu-2$ at~$t=t_0$.
Thus in space dimension~$n=1$ Proposition~\ref{P-Na} can be precised.

\begin{Proposition} \label{P-Na1}
In space dimension $n=1$,
\eqref{E-NDD} is equivalent to~\eqref{E-Ncl}
and can be written in the following form:
there exist~$t_1,\dotsc,t_\nu$ such that
\begin{equation} \label{E-NDD1}
\prod_{j=1}^\nu (t-t_j)^2 \, \Bigl|\widecheck{N}\bigr(t,\sigma_1(t,\xi),\xi\bigr)\Bigr|
  \lesssim  \bigl|\sigma_2(t,\xi)-\sigma_1(t,\xi)\bigr| \, .
\end{equation}
\end{Proposition}

\begin{proof}
We prove that if $n=1$,
\eqref{E-Ncpp} implies~\eqref{E-NDD1}.
As before, with no loss of~generality,
we can assume that $\Delta_{\partial_t L}$ vanishes only in~$0$,
and there exists $r \ge 1$ such that
\[
\bigl|\sigma_2(t)-\sigma_1(t)\bigl| \gtrsim t^r \, .
\]
Hence, for $t\ge\e^{1/r}$, $\e=\frac{1}{|\xi|}$ and $|\xi|\ge1$
we have
\[
\bigl|\sigma_2(t)-\sigma_1(t)\bigl|\,|\xi| \gtrsim 1 \, ,
\]
and, consequently,
\[
\frac{\bigl|\sigma_2(t)-\sigma_1(t)\bigl|\,|\xi|}
     {\bigl|\sigma_2(t)-\sigma_1(t)\bigl|\,|\xi| + 1} \gtrsim 1 \, .
\]

If \eqref{E-NDD1} fails to hold then there exists $m\ge 3$
such that
\[
\frac{\widecheck{N}\bigl(t,\sigma_j(t,\xi),\xi\bigr)}
     {\sigma_2(t,\xi)-\sigma_1(t,\xi)}
  \approx  \frac{1}{t^m}
\quad  \text{for $j=1$ or $j=2$} \, ,
\]
and we have
\[
\int_{\e^{1/r}}^T
    \sqrt{\frac{\bigl|\widecheck{N}\bigl(t,\sigma_j(t,\xi),\xi\bigr)\bigr|}
               {\bigl|\sigma_3(t,\xi)-\sigma_1(t,\xi)\bigr|+1} \,} \, dt
  \gtrsim \int_{\e^{1/r}}^T \frac{1}{t^{m/2}} \, dt
  \approx |\xi|^{\frac{m/2-1}{r}} \, ,
\]
thus~\eqref{E-Ncpp} cannot be satisfied.
\end{proof}

\subsection*{Case~II: $\Delta_L\equiv 0$ and $\Delta_L^{(1)}\protect\not\equiv 0$}
\text{}

With no loss of~generality,
we can assume that $\tau_1\equiv\tau_2$ and $\tau_3\not\equiv\tau_1$.

\begin{Proposition} \label{P-3}
Assume that $\tau_1\equiv\tau_2$ and $\tau_3\not\equiv\tau_1$.
If
\begin{subequations}
\begin{equation} \label{E-80}
\widecheck{M}\bigr(t,\tau_1(t,\xi),\xi\bigr)
  \equiv  0
\end{equation}
and
\begin{equation} \label{E-QD}
\Bigl|Q\bigr(t,\tau_j(t,\xi),\xi\bigr)\Bigr|
  \lesssim  \sqrt{\Delta_L^{(1)}(t,\xi) \,}
            +  \Bigl|\partial_t \sqrt{\Delta_L^{(1)}(t,\xi) \,}\Bigr| \, ,
\end{equation}
\end{subequations}
for $j=1,3$,
where $Q(t,\tau,\xi)= \frac{\widecheck{M}(t,\tau,\xi)}{\tau - \tau_1(t,\xi)}$,
then condition~\eqref{E-Mcl} is verified.
\end{Proposition}

\begin{proof}
The proof is similar to that of~Proposition~\ref{P-Ma}:
we prove that~\eqref{E-80} and~\eqref{E-QD} implies~\eqref{E-Mct} and~\eqref{E-dMct},
hence, by~Proposition~\ref{P-M} we get \eqref{E-Mcl}.

Note that, by~\eqref{E-80},
\eqref{E-Mct} with $j=1$ or $j=2$ is trivially satisfied,
thus we need only to prove~\eqref{E-Mct} with $j=3$.

Let $\e=|\xi|^{-2}$
with $|\xi|$ large enough,
and let $A_{\xi,\e}$ be the set given by Lemma~\ref{L-Prato 2}
with $\Delta(t,\xi) = \Delta_L^{(1)}(t,\xi)$;
we~have
\[
\int_{A_{\xi,\e}} \frac{\bigl|\widecheck{M}(\tau_3)\bigr|}
                  {\bigl(|\tau_3 - \tau_1|+1\bigr)^2} \, dt
  \lesssim  \int_{A_{\xi,\e}} |\xi|^2 \, dt
  \lesssim  \e\,|\xi|^2
   = 1 \, .
\]

On the other side,
as
\[
\frac{\bigl|\widecheck{M}(\tau_3)\bigr|}
     {\bigl(|\tau_3 - \tau_1|+1\bigr)^2}
  \le  \frac{\bigl|Q(\tau_3)\bigr|}
            {\sqrt{\Delta_L^{(1)} \,}} \, ,
\]
assuming~\eqref{E-QD},
thanks to Lemma~\ref{L-Prato 2},
we have
\[
\int_{\left]-\delta,T+\delta\right[ \setminus A_{\xi,\e}}
             \frac{\bigl|\widecheck{M}(\tau_3)\bigr|}
                  {\bigl(|\tau_3-\tau_1|+1\bigr)^2} \, dt
  \lesssim  \int_{\left]-\delta,T+\delta\right[ \setminus A_{\xi,\e}}
                  1 + \frac{\bigl|\partial_t \Delta_L^{(1)}(t,\xi)\bigr|}
                           {\Delta_L^{(1)}(t,\xi)} \, dt
  \logestimate \, .
\]

Now we prove that~\eqref{E-dMct} holds true.

As $\tau_3\not\equiv\tau_1$,
by Lagrange interpolation formula,
we~get
\begin{align}
Q(\tau)
&  =  \frac{Q(\tau_3)}{\tau_3-\tau_1} \, (\tau-\tau_1)
      +  \frac{Q(\tau_1)}{\tau_1-\tau_3} \, (\tau-\tau_3) \notag
\intertext{%
for a.e.~$(t,\xi)$,
hence}
\widecheck{M}(\tau)
&  =  \frac{Q(\tau_3)}{\tau_3-\tau_1} \, L_{12}
      +  \frac{Q(\tau_1)}{\tau_1-\tau_3} \, L_{13} \, , \label{E-555}
\intertext{%
and, consequently,}
\partial_\tau \widecheck{M}(\tau)
&  =  \frac{Q(\tau_3)}{\tau_3-\tau_1} \, (L_1+L_2)
      +  \frac{Q(\tau_1)}{\tau_1-\tau_3} \, (L_1+L_3) \, . \label{E-678}
\end{align}

Let $\e=|\xi|^{-1}$,
with $|\xi|$ large enough,
and let $A_{\xi,\e}$ be as above,
we have
\begin{align*}
\int_{A_{\xi,\e}} \frac{\bigl|\partial_\tau \widecheck{M}(\tau_1)\bigr|
                         +  \bigl|\partial_\tau \widecheck{M}(\tau_3)\bigr|}
                       {|\tau_1-\tau_3|+1} \, dt
&  \lesssim  \int_{A_{\xi,\e}} |\xi| \, dt
   \lesssim  \e\,|\xi|
   = 1 \, ,
\intertext{%
whereas, thanks to~\eqref{E-678} and Lemma~\ref{L-Prato 2} we have}
\int_{\left]-\delta,T+\delta\right[ \setminus A_{\xi,\e}}
             \frac{\bigl|\partial_\tau \widecheck{M}(\tau_1)\bigr|
                    +  \bigl|\partial_\tau \widecheck{M}(\tau_3)\bigr|}
                  {|\tau_1-\tau_3|+1} \, dt
&  \lesssim  \int_{\left]-\delta,T+\delta\right[ \setminus A_{\xi,\e}}
             \frac{\bigl|Q(\tau_1)\bigr| + \bigl|Q(\tau_3)\bigr|}
                  {|\tau_1-\tau_3|+1} \, dt
\intertext{%
hence, by~\eqref{E-QD},}
\int_{\left]-\delta,T+\delta\right[ \setminus A_{\xi,\e}}
             \frac{\bigl|\partial_\tau \widecheck{M}(\tau_1)\bigr|
                    +  \bigl|\partial_\tau \widecheck{M}(\tau_3)\bigr|}
                  {|\tau_1-\tau_3|+1} \, dt
&  \lesssim  \int_{\left]-\delta,T+\delta\right[ \setminus A_{\xi,\e}}
                 1 + \frac{\bigl|\partial_\tau \Delta_L^{(1)}(t,\xi)\bigr|}
                          {\Delta_L^{(1)}(t,\xi)} \, dt  \\
&  \logestimate \, . \qedhere
\end{align*}
\end{proof}

\begin{Proposition}
In space dimension $n=1$,
\eqref{E-Mcl} is equivalent to~\eqref{E-80} and~\eqref{E-QD}.

Moreover~\eqref{E-QD} can be written in the following form:
there exist~$t_1,\dotsc,t_\nu$ such that
\begin{equation} \label{E-MDD1}
\prod_{j=1}^\nu |t-t_j| \, \Bigl|Q\bigr(t,\tau_k(t,\xi),\xi\bigr)\Bigr|
  \lesssim  \bigl|\tau_3(t,\xi)-\tau_1(t,\xi)\bigr| \, ,
\qquad
k=1,3 \, .
\end{equation}
\end{Proposition}

\begin{Remark}
In the special case $\nu=1$ and $t_1=0$,
condition~\eqref{E-MDD1}
reduces to
\begin{equation} \label{E-81}
|t| \, \Bigl|Q\bigr(t,\tau_k(t,\xi),\xi\bigr)\Bigr|
  \lesssim  \bigl|\tau_3(t,\xi)-\tau_1(t,\xi)\bigr| \, ,
\qquad
k=1,3 \, .
\end{equation}
\end{Remark}

\begin{proof}
We prove that if $\tau_1\equiv\tau_2$,
then~\eqref{E-80} is necessary in order to have~\eqref{E-Mcl}.

Indeed if \eqref{E-80} fails to hold then
$\widecheck{M}\bigr(t,\tau_1(t,\xi),\xi\bigr) \approx |\xi|^2$
and we~have
\[
\frac{\widecheck{M}(\tau_1)}
     {\bigl(|\tau_1-\tau_2| + 1\bigr)
      \bigl(|\tau_1-\tau_3| + 1\bigr)}
  =  \frac{\widecheck{M}(\tau_1)}
          {|\tau_1-\tau_3| + 1}
  \approx  \frac{|\xi|^2}{|\xi| + 1} \, ,
\]
hence~\eqref{E-Mct} cannot be verified.

As before, we can assume that $\Delta_L^{(1)}$ vanishes
only in~$0$,
so that~\eqref{E-MDD1} reduces to~\eqref{E-81}.
Moreover we can assume that there exists $r \ge 1$ such that
for $t\ge\e^{1/r}$, $\e=\frac{1}{|\xi|}$ and $|\xi|\ge1$
we have
\[
\bigl|\tau_3(t)-\tau_1(t)\bigl|\,|\xi| \gtrsim 1 \, .
\]

If \eqref{E-81} fails to hold then there exists $m\ge2$
such that
\[
\frac{\widecheck{Q}\bigl(t,\tau_3(t,\xi),\xi\bigr)}
     {\tau_3(t,\xi)-\tau_1(t,\xi)}
  \approx  \frac{1}{t^m} \, ,
\]
and, as
\[
\widecheck{M}(t,\tau,\xi)
  =  Q(t,\tau,\xi) \, \bigl( \tau - \tau_1(t,\xi) \bigr)
\]
we have
\[
\partial_\tau \widecheck{M}(t,\tau_1,\xi)
  =  Q(t,\tau_1,\xi) \, ,
\]
hence
\[
\int_{\e^{1/r}}^T
    \frac{\bigl|\partial_\tau \widecheck{M}(\tau_1)\bigr|}
         {|\tau_3-\tau_1|+1} \, dt
  =  \int_{\e^{1/r}}^T
       \frac{\bigl|Q(\tau_1)\bigr|}
            {|\tau_3-\tau_1|+1} \, dt
  \gtrsim \int_{\e^{1/r}}^T \frac{1}{t^m} \, dt
  \approx |\xi|^{\frac{m-1}{r}} \, .
\]
Thus~\eqref{E-dMct} cannot be satisfied.
\end{proof}

Note that we need not to assume $n=1$
to prove the necessity of~\eqref{E-80}.

\medskip

The term of~order 1
is treated as in the previous case.

\subsection*{Case~III: $\Delta_L\equiv\Delta_L^{(1)}\equiv0$
   (operators with triple characteristics of~constant multiplicity)}
\text{}

\begin{Proposition} \label{P-4}
If $\Delta_L\equiv\Delta_L^{(1)}\equiv0$
then $L$ has a unique triple root:
\begin{equation} \label{E-L3}
L(t,\tau,\xi)
  =  \bigl(\tau-\tau_1(t,\xi)\bigr)^3 \, ,
\end{equation}
and Hypothesis~\eqref{E-1} and~\eqref{E-2} are satisfied.
Moreover Hypothesis~\eqref{E-3}, \eqref{E-4}, \eqref{E-Mcl} and \eqref{E-Ncl}
are satisfied if, and only if,
\begin{equation} \label{E-charconst}
\widecheck{M}(\tau_1) \equiv 0 \, , \quad
\partial_\tau \widecheck{M}(\tau_1) \equiv 0 \, , \quad
\widecheck{N}(\tau_1) \equiv 0 \, .
\end{equation}
\end{Proposition}

\begin{proof}
It's clear that if $\Delta_L\equiv\Delta_L^{(1)}\equiv0$
then $L$ reduces to~\eqref{E-L3}.

If $L$ is as in~\eqref{E-L3}
we have
\[
\mathcal{L}(t,\tau,\xi)
   =  \bigl(\tau-\tau_1(t,\xi)\bigr)^3 - 6\,\bigl(\tau-\tau_1(t,\xi)\bigr)
\]
and
\[
\lambda_1(t,\xi) = \tau_1(t,\xi) \, , \qquad
\lambda_2(t,\xi) = \tau_1(t,\xi)+\sqrt{6\,} \, , \qquad
\lambda_3(t,\xi) = \tau_1(t,\xi)-\sqrt{6\,} \, ,
\]
and it's clear that \eqref{E-1} and~\eqref{E-2} are satisfied.

It's also clear that~\eqref{E-charconst} holds true
if, and only if,
\[
\widecheck{M}(t,\tau,\xi)
   =  m_0(t) \, \bigl(\tau-\tau_1(t,\xi)\bigr)^2
\quad\text{and}\quad
\widecheck{N}(t,\tau,\xi)
   =  n_0(t) \, \bigl(\tau-\tau_1(t,\xi)\bigr) \, .
\]
Thus~\eqref{E-3}, \eqref{E-4}, \eqref{E-Mcl} and~\eqref{E-Ncl} hold true.

On the converse,
if $\widecheck{M}(\tau_1)\not\equiv0$ then
\eqref{E-Mct} is not satisfied since
\[
\int_0^T \bigl|\widecheck{M}\bigl(t,\tau_1(t,\xi),\xi\bigr)\bigr| \, dt
  \approx  |\xi|^2 \, .
\]

Analogously,
if $\partial_\tau \widecheck{M}(\tau_1)\not\equiv0$ then
\eqref{E-dMct} is not satisfied since
\[
\int_0^T \Bigl|\partial_\tau \widecheck{M}\bigl(t,\tau_1(t,\xi),\xi\bigr)\Bigr| \, dt
  \approx  |\xi| \, .
\]

Finally,
if $\widecheck{N}(\mu_1)\not\equiv0$ then
\eqref{E-Ncpp} is not satisfied since
\[
\int_0^T \sqrt{ \Bigl|\widecheck{N}\bigl(t,\tau_1(t,\xi),\xi\bigr)\Bigr|\,} \, dt
  \approx  |\xi|^{1/2} \, . \qedhere
\]
\end{proof}

Conditions~\eqref{E-charconst}
correspond to the conditions of~\emph{good decomposition}~\cite{DeParis}:
the operator $P$ can be written as
\[
P
  =  L_1^3 + m_0(t) \, L_1^2 + n_0(t) \, L_1 + p_0 \, ,
\]
where $L_1^3 = L_1\circ L_1\circ L_1$
and $L_1^2 = L_1\circ L_1$.

To see this, following~\cite{TV2008},
we have to check two conditions.
The first is
\begin{align}
&  \text{the principal symbol of~\ $P - L_1^3$ \
         is divisible by \ $m_0(\tau-\tau_1)^2$} \, , \label{L-Levi1}
\intertext{%
and the second condition is}
&  \text{the principal symbol of~\ $P - (L_1^3 + m_0(t) \, L_1^2)$ \
         is divisible by \ $(\tau-\tau_1)$} \, . \label{L-Levi2}
\end{align}

As $\tau_1\equiv\tau_2\equiv\tau_3$,
the operator~$\wt{L}_{123,\e}$ in~\eqref{E-wtLLL}
reduces to~$L_1^3$ if~$\e=0$.
According to~\eqref{E-dec} we have
\[
P
  =  L_1^3 + \widecheck{M} + \frac{1}{2} \partial_t\partial_\tau \widecheck{M} + \widecheck{N} + p \, ,
\]
thus~\eqref{L-Levi1} holds true if, and only if,
the first two conditions in~\eqref{E-charconst} hold true.
In~this case
\[
\widecheck{M}(t,\tau,\xi)
  =  m_0(t) \, L_{1,1}(t,\tau,\xi)
  =  m_0(t) \, \bigl(\tau - \tau_1(t,\xi)\bigr)^2 \, .
\]

As $\tau_1\equiv\tau_2$,
the operator~$\wt{L}_{12,\e}$ in~\eqref{E-wtLL}
reduces to~$L_1^2$ if~$\e=0$.
By~\eqref{E-tLjh} with $\e=0$ we have
\begin{align*}
L_1^2
&  =  L_{11}
       + \frac{1}{2} \, \partial_t \partial_\tau L_{11} \, ,
\intertext{and,
multiplying by $m_0$, we get}
m_0\,L_1^2
&  =  \widecheck{M}
      + \frac{1}{2} \, \partial_t \partial_\tau \widecheck{M}
      - m_0'(t) \, L_1 \, ,
\end{align*}
hence
\[
P
  =  L_1^3 + m_0\,L_1^2
     +  m_0'(t) \, L_1 + \widecheck{N} + p \, ,
\]
thus~\eqref{L-Levi2} holds true if, and only if,
the third condition in~\eqref{E-charconst} holds true.

\section{Operators with constant coefficients principal part} \label{S-cc}

In this section we prove the following Proposition

\begin{Proposition}
Assume that the coefficients are analytic
and those of~the principal symbol are constant.

Then Hypothesis~\eqref{E-1}, \eqref{E-2}, \eqref{E-3} and~\eqref{E-4} are satisfied,
whereas~\eqref{E-Mcl} and~\eqref{E-Ncl} are necessary and sufficient
for the $\mathcal{C}^\infty$ well-posedness.
\end{Proposition}

\begin{proof}
If the coefficients are analytic then,
Hypothesis~\eqref{E-1}, \eqref{E-2}, \eqref{E-3} and~\eqref{E-4} are satisfied
(cf.~\textsection\ref{S-an}).

\smallskip

Now, we recall that if the coefficients are constant,
the~necessary and sufficient conditions for the $\mathscr{C}^\infty$
well-posedness is well-known,
see~\cite{Garding} \cite{Hormander} \cite{Svensson}:
\begin{equation} \label{E-Garding}
\text{there exists~$C>0$ such that
      $\tau^3 + \sum_{j+|\alpha| \le 3} a_{j,\alpha}(t)\tau^j\xi^\alpha\ne0$ 
      if~$\xi\in\R$ and $|\Im\tau|>C$} \, .
\end{equation}

Various equivalent conditions have been given.
According to~\cite{Svensson},
\eqref{E-Garding} is equivalent to the following conditions (cf.~\cite{Peyser1963}):

\smallskip

there exist bounded functions $\ell_1,\ell_2,\ell_3,\mathsf{m}_1,\mathsf{m}_2$
such that
\begin{align}
M(t,\tau,\xi)
&  =  \ell_1(t,\xi) \, \bigl(\tau - \tau_2(t,\xi)\bigr) \bigl(\tau - \tau_3(t,\xi)\bigr)
      +  \ell_2(t,\xi) \, \bigl(\tau - \tau_3(t,\xi)\bigr) \bigl(\tau - \tau_1(t,\xi)\bigr) \label{E-849}  \\
&  \qquad  +  \ell_3(t,\xi) \, \bigl(\tau - \tau_1(t,\xi)\bigr) \bigl(\tau - \tau_2(t,\xi)\bigr) \notag  \\
N(t,\tau,\xi)
&  =  \mathsf{m}_1(t,\xi) \, \bigl(\tau - \sigma_2(t,\xi)\bigr)
      +  \mathsf{m}_2(t,\xi) \, \bigl(\tau - \sigma_1(t,\xi)\bigr) \label{E-847}
\end{align}
for a.e.~$(t,\xi)$.

We recall that the above conditions are also necessary and sufficient
for the $\mathcal{C}^\infty$ well-posedness
if the coefficients of~the lower order terms can vary~\cite{Dunn}, \cite{Wakabayashi1980}
(also if the coefficients of the lower order terms are only~$\mathcal{C}^\infty$).

As we have proved in~\textsection\ref{S-proof}
that our conditions are sufficient for the well-posedness,
whereas conditions~\eqref{E-849} and~\eqref{E-847}
are necessary,
it remains to prove that
if~\eqref{E-849} and~\eqref{E-847} holds true
then~\eqref{E-Mcl} and~\eqref{E-Ncl} are satisfied.

\smallskip

We prove at first that~\eqref{E-849} implies~\eqref{E-Mcl}.
We have to distinguish three cases as in the previous section.

If $\Delta_L\not\equiv0$,
then~\eqref{E-849} is~\eqref{E-845}
and the $\ell_j$ are given by~\eqref{E-ell}.
The boundness of~the $\ell_j$ is equivalent to~\eqref{E-170D},
hence, by~Proposition~\ref{P-Ma}, \eqref{E-Mcl} holds true.

If $\Delta_L\equiv0$ and $\Delta_{\partial_\tau}L\not\equiv0$,
then 
with no loss of~generality,
we can assume that $\tau_1\equiv\tau_2$ and $\tau_3\not\equiv\tau_1$.
In this case the right hand side of~\eqref{E-849}
is divisible by $\tau-\tau_1$,
hence~$M$ must verify~\eqref{E-80}.
Moreover~\eqref{E-849} reduces to~\eqref{E-555},
thus the boundness of~the $\ell_j$ is equivalent to~\eqref{E-QD}.
By~Proposition~\ref{P-3} we get~\eqref{E-Mcl}.

If $\Delta_{\partial_\tau}L\equiv0$,
then $\tau_1\equiv\tau_2\equiv\tau_3$ and the right hand side of~\eqref{E-849}
is divisible by $(\tau-\tau_1)^2$,
hence~$M$ must verify the first two conditions in~\eqref{E-charconst},
and, by~Proposition~\ref{P-4}, we get~\eqref{E-Mcl}.

\smallskip

Now we prove that~\eqref{E-847} implies~\eqref{E-Ncl}.

{}From~\eqref{E-847} we have
\[
N(\sigma_1)
  =  \mathsf{m}_2 \, (\sigma_1 - \sigma_2) \, ,
\qquad
N(\sigma_2)
  =  \mathsf{m}_1 \, (\sigma_2 - \sigma_1) \, ,
\]
from which we deduce that
if $\Delta_{\partial_\tau L}\not\equiv0$ then~\eqref{E-NDD} is satisfied;
by~Proposition~\ref{P-Na}, \eqref{E-Ncl} holds true.
On~the other side, if $\Delta_{\partial_\tau L}\equiv0$
then the third condition in~\eqref{E-charconst} holds true,
hence, by~Proposition~\ref{P-4}, we get~\eqref{E-Mcl}.
\end{proof}

\appendix
\section{}

\begin{Lemma} \label{L-A-1}
Let
\[
p(\tau)
  =  \tau^3 + A_1 \tau^2 + A_2 \tau + A_3
  =  (\tau-\tau_1) \, (\tau-\tau_2) \, (\tau-\tau_3) \, ,
\]
and
\[
p'(\tau)
  =  3\,\tau^2 + 2 A_1 \tau + A_2
  =  3 \, (\tau-\sigma_1) \, (\tau-\sigma_2) \, ,
\]
then
\begin{equation} \label{E-A1}
(\tau_{1}-\tau_{2})^2+(\tau_{2}-\tau_{3})^2+(\tau_{3}-\tau_{1})^2
  =  \frac{9}{2} \, (\sigma_{2}-\sigma_{1})^2 \, .
\end{equation}
\end{Lemma}

\begin{proof}
Using Vieta's formulas
\begin{align*}
(\tau_{1}-\tau_{2})^2 +
& (\tau_{2}-\tau_{3})^2+(\tau_{3}-\tau_{1})^2  \\
&  =  2\,\bigl[\tau_{1}^2+\tau_{2}^2+\tau_{3}^2
               -\tau_{1}\tau_{2}-\tau_{2}\tau_{3}-\tau_{3}\tau_{1}\bigr]  \\
&  =  2\,\bigl[(\tau_{1}+\tau_{2}+\tau_{3})^2
               -3(\tau_{1}\tau_{2}+\tau_{2}\tau_{3}+\tau_{3}\tau_{1})\bigr]  \\
&  =  2\,[A_1^2 - 3A_2] \, ,
\end{align*}
whereas
\[
(\sigma_{2}-\sigma_{1})^2
  =  \frac{4}{9}\,A_1^2 - \frac{4}{3}\,A_2 \, ,
\]
from which we get the result.
\end{proof}

\providecommand{\bysame}{\leavevmode\hbox to3em{\hrulefill}\thinspace}
\providecommand{\MR}{\relax\ifhmode\unskip\space\fi MR }
\providecommand{\MRhref}[2]{%
  \href{http://www.ams.org/mathscinet-getitem?mr=#1}{#2}
}
\providecommand{\href}[2]{#2}

\end{document}